%%%%%%%%%%%%%%%%%%%%%%%%%%%%%%%%%%%%%%%%%%%%%%%%%%%%%%%%%%%%%%%%%%%%%%%%%%%
%%% Title:  Product preserving functors of infinite dimensional manifolds
%%% Author: Andreas Kriegl, Peter W. Michor
%%% Remark: AmSTeX and diag.tex
%%% Series: 
%%%%%%%%%%%%%%%%%%%%%%%%%%%%%%%%%%%%%%%%%%%%%%%%%%%%%%%%%%%%%%%%%%%%%%%%%%%
% TeX-NMB program applied to the file from 1996.11.4;14:20 on 1996.11.4; 14:23
% TeX-NMB program applied to the file from 1996.10.10;12:50 on 1996.10.10; 12:51
% TeX-NMB program applied to the file from 1996.10.4;7:42 on 1996.10.4; 7:42
% TeX-NMB program applied to the file from 1996.9.20;16:41 on 1996.9.20; 16:41
% TeX-NMB program applied to the file from 1996.9.20;16:25 on 1996.9.20; 16:41
\input amstex
\input amsppt.sty   
%\input diag.tex
% macro package for drawing diagrams --
% essentially a part of lamstex.tex
% use exactly following the lamstex manual 
% with '..CD..' replaced by  '..newCD..'
% the original construction $\CD ... \endCD$ remains without any 
% change
%
\catcode`\@=11

\def\East#1#2{\setboxz@h{$\m@th\ssize\;{#1}\;\;$}%
 \setbox@ne\hbox{$\m@th\ssize\;{#2}\;\;$}\setbox\tw@\hbox{$\m@th#2$}%
 \dimen@\minaw@
 \ifdim\wdz@>\dimen@ \dimen@\wdz@ \fi  \ifdim\wd@ne>\dimen@ \dimen@\wd@ne \fi
 \ifdim\wd\tw@>\z@
  \mathrel{\mathop{\hbox to\dimen@{\rightarrowfill}}\limits^{#1}_{#2}}%
 \else
  \mathrel{\mathop{\hbox to\dimen@{\rightarrowfill}}\limits^{#1}}%
 \fi}
\def\West#1#2{\setboxz@h{$\m@th\ssize\;\;{#1}\;$}%
 \setbox@ne\hbox{$\m@th\ssize\;\;{#2}\;$}\setbox\tw@\hbox{$\m@th#2$}%
 \dimen@\minaw@
 \ifdim\wdz@>\dimen@ \dimen@\wdz@ \fi \ifdim\wd@ne>\dimen@ \dimen@\wd@ne \fi
 \ifdim\wd\tw@>\z@
  \mathrel{\mathop{\hbox to\dimen@{\leftarrowfill}}\limits^{#1}_{#2}}%
 \else
  \mathrel{\mathop{\hbox to\dimen@{\leftarrowfill}}\limits^{#1}}%
 \fi}
\font\arrow@i=lams1
\font\arrow@ii=lams2
\font\arrow@iii=lams3
\font\arrow@iv=lams4
\font\arrow@v=lams5
\newbox\zer@
\newdimen\standardcgap
\standardcgap=40\p@
\newdimen\hunit
\hunit=\tw@\p@
\newdimen\standardrgap
\standardrgap=32\p@
\newdimen\vunit
\vunit=1.6\p@
\def\Cgaps#1{\RIfM@
  \standardcgap=#1\standardcgap\relax \hunit=#1\hunit\relax
 \else \nonmatherr@\Cgaps \fi}
\def\Rgaps#1{\RIfM@
  \standardrgap=#1\standardrgap\relax \vunit=#1\vunit\relax
 \else \nonmatherr@\Rgaps \fi}
\newdimen\getdim@
\def\getcgap@#1{\ifcase#1\or\getdim@\z@\else\getdim@\standardcgap\fi}
\def\getrgap@#1{\ifcase#1\getdim@\z@\else\getdim@\standardrgap\fi}
\def\cgaps#1{\RIfM@
 \cgaps@{#1}\edef\getcgap@##1{\i@=##1\relax\the\toks@}\toks@{}\else
 \nonmatherr@\cgaps\fi}
\def\rgaps#1{\RIfM@
 \rgaps@{#1}\edef\getrgap@##1{\i@=##1\relax\the\toks@}\toks@{}\else
 \nonmatherr@\rgaps\fi}
\def\Gaps@@{\gaps@@}
\def\cgaps@#1{\toks@{\ifcase\i@\or\getdim@=\z@}%
 \gaps@@\standardcgap#1;\gaps@@\gaps@@
 \edef\next@{\the\toks@\noexpand\else\noexpand\getdim@\noexpand\standardcgap
  \noexpand\fi}%
 \toks@=\expandafter{\next@}}
\def\rgaps@#1{\toks@{\ifcase\i@\getdim@=\z@}%
 \gaps@@\standardrgap#1;\gaps@@\gaps@@
 \edef\next@{\the\toks@\noexpand\else\noexpand\getdim@\noexpand\standardrgap
  \noexpand\fi}%
 \toks@=\expandafter{\next@}}
\def\gaps@@#1#2;#3{\mgaps@#1#2\mgaps@
 \edef\next@{\the\toks@\noexpand\or\noexpand\getdim@
  \noexpand#1\the\mgapstoks@@}%
 \global\toks@=\expandafter{\next@}%
 \DN@{#3}%
 \ifx\next@\Gaps@@\gdef\next@##1\gaps@@{}\else
  \gdef\next@{\gaps@@#1#3}\fi\next@}
\def\mgaps@#1{\let\mgapsnext@#1\FN@\mgaps@@}
\def\mgaps@@{\ifx\next\space@\DN@. {\FN@\mgaps@@}\else
 \DN@.{\FN@\mgaps@@@}\fi\next@.}
\def\mgaps@@@{\ifx\next\w\let\next@\mgaps@@@@\else
 \let\next@\mgaps@@@@@\fi\next@}
\newtoks\mgapstoks@@
\def\mgaps@@@@@#1\mgaps@{\getdim@\mgapsnext@\getdim@#1\getdim@
 \edef\next@{\noexpand\getdim@\the\getdim@}%
 \mgapstoks@@=\expandafter{\next@}}
\def\mgaps@@@@\w#1#2\mgaps@{\mgaps@@@@@#2\mgaps@
 \setbox\zer@\hbox{$\m@th\hskip15\p@\tsize@#1$}%
 \dimen@\wd\zer@
 \ifdim\dimen@>\getdim@ \getdim@\dimen@ \fi
 \edef\next@{\noexpand\getdim@\the\getdim@}%
 \mgapstoks@@=\expandafter{\next@}}
\def\changewidth#1#2{\setbox\zer@\hbox{$\m@th#2}%
 \hbox to\wd\zer@{\hss$\m@th#1$\hss}}
\atdef@({\FN@\ARROW@}
\def\ARROW@{\ifx\next)\let\next@\OPTIONS@\else
 \DN@{\csname\string @(\endcsname}\fi\next@}
\newif\ifoptions@
\def\OPTIONS@){\ifoptions@\let\next@\relax\else
 \DN@{\options@true\begingroup\optioncodes@}\fi\next@}
\newif\ifN@
\newif\ifE@
\newif\ifNESW@
\newif\ifH@
\newif\ifV@
\newif\ifHshort@
\expandafter\def\csname\string @(\endcsname #1,#2){%
 \ifoptions@\let\next@\endgroup\else\let\next@\relax\fi\next@
 \N@false\E@false\H@false\V@false\Hshort@false
 \ifnum#1>\z@\E@true\fi
 \ifnum#1=\z@\V@true\tX@false\tY@false\a@false\fi
 \ifnum#2>\z@\N@true\fi
 \ifnum#2=\z@\H@true\tX@false\tY@false\a@false\ifshort@\Hshort@true\fi\fi
 \NESW@false
 \ifN@\ifE@\NESW@true\fi\else\ifE@\else\NESW@true\fi\fi
 \arrow@{#1}{#2}%
 \global\options@false
 \global\scount@\z@\global\tcount@\z@\global\arrcount@\z@
 \global\s@false\global\sxdimen@\z@\global\sydimen@\z@
 \global\tX@false\global\tXdimen@i\z@\global\tXdimen@ii\z@
 \global\tY@false\global\tYdimen@i\z@\global\tYdimen@ii\z@
 \global\a@false\global\exacount@\z@
 \global\x@false\global\xdimen@\z@
 \global\X@false\global\Xdimen@\z@
 \global\y@false\global\ydimen@\z@
 \global\Y@false\global\Ydimen@\z@
 \global\p@false\global\pdimen@\z@
 \global\label@ifalse\global\label@iifalse
 \global\dl@ifalse\global\ldimen@i\z@
 \global\dl@iifalse\global\ldimen@ii\z@
 \global\short@false\global\unshort@false}
\newif\iflabel@i
\newif\iflabel@ii
\newcount\scount@
\newcount\tcount@
\newcount\arrcount@
\newif\ifs@
\newdimen\sxdimen@
\newdimen\sydimen@
\newif\iftX@
\newdimen\tXdimen@i
\newdimen\tXdimen@ii
\newif\iftY@
\newdimen\tYdimen@i
\newdimen\tYdimen@ii
\newif\ifa@
\newcount\exacount@
\newif\ifx@
\newdimen\xdimen@
\newif\ifX@
\newdimen\Xdimen@
\newif\ify@
\newdimen\ydimen@
\newif\ifY@
\newdimen\Ydimen@
\newif\ifp@
\newdimen\pdimen@
\newif\ifdl@i
\newif\ifdl@ii
\newdimen\ldimen@i
\newdimen\ldimen@ii
\newif\ifshort@
\newif\ifunshort@
\def\zero@#1{\ifnum\scount@=\z@
 \if#1e\global\scount@\m@ne\else
 \if#1t\global\scount@\tw@\else
 \if#1h\global\scount@\thr@@\else
 \if#1'\global\scount@6 \else
 \if#1`\global\scount@7 \else
 \if#1(\global\scount@8 \else
 \if#1)\global\scount@9 \else
 \if#1s\global\scount@12 \else
 \if#1H\global\scount@13 \else
 \Err@{\Invalid@@ option \string\0}\fi\fi\fi\fi\fi\fi\fi\fi\fi
 \fi}
\def\one@#1{\ifnum\tcount@=\z@
 \if#1e\global\tcount@\m@ne\else
 \if#1h\global\tcount@\tw@\else
 \if#1t\global\tcount@\thr@@\else
 \if#1'\global\tcount@4 \else
 \if#1`\global\tcount@5 \else
 \if#1(\global\tcount@10 \else
 \if#1)\global\tcount@11 \else
 \if#1s\global\tcount@12 \else
 \if#1H\global\tcount@13 \else
 \Err@{\Invalid@@ option \string\1}\fi\fi\fi\fi\fi\fi\fi\fi\fi
 \fi}
\def\a@#1{\ifnum\arrcount@=\z@
 \if#10\global\arrcount@\m@ne\else
 \if#1+\global\arrcount@\@ne\else
 \if#1-\global\arrcount@\tw@\else
 \if#1=\global\arrcount@\thr@@\else
 \Err@{\Invalid@@ option \string\a}\fi\fi\fi\fi
 \fi}
\def\ds@(#1;#2){\ifs@\else
 \global\s@true
 \sxdimen@\hunit \global\sxdimen@#1\sxdimen@\relax
 \sydimen@\vunit \global\sydimen@#2\sydimen@\relax
 \fi}
\def\dtX@(#1;#2){\iftX@\else
 \global\tX@true
 \tXdimen@i\hunit \global\tXdimen@i#1\tXdimen@i\relax
 \tXdimen@ii\vunit \global\tXdimen@ii#2\tXdimen@ii\relax
 \fi}
\def\dtY@(#1;#2){\iftY@\else
 \global\tY@true
 \tYdimen@i\hunit \global\tYdimen@i#1\tYdimen@i\relax
 \tYdimen@ii\vunit \global\tYdimen@ii#2\tYdimen@ii\relax
 \fi}
\def\da@#1{\ifa@\else\global\a@true\global\exacount@#1\relax\fi}
\def\dx@#1{\ifx@\else
 \global\x@true
 \xdimen@\hunit \global\xdimen@#1\xdimen@\relax
 \fi}
\def\dX@#1{\ifX@\else
 \global\X@true
 \Xdimen@\hunit \global\Xdimen@#1\Xdimen@\relax
 \fi}
\def\dy@#1{\ify@\else
 \global\y@true
 \ydimen@\vunit \global\ydimen@#1\ydimen@\relax
 \fi}
\def\dY@#1{\ifY@\else
 \global\Y@true
 \Ydimen@\vunit \global\Ydimen@#1\Ydimen@\relax
 \fi}
\def\p@@#1{\ifp@\else
 \global\p@true
 \pdimen@\hunit \divide\pdimen@\tw@ \global\pdimen@#1\pdimen@\relax
 \fi}
\def\L@#1{\iflabel@i\else
 \global\label@itrue \gdef\label@i{#1}%
 \fi}
\def\l@#1{\iflabel@ii\else
 \global\label@iitrue \gdef\label@ii{#1}%
 \fi}
\def\dL@#1{\ifdl@i\else
 \global\dl@itrue \ldimen@i\hunit \global\ldimen@i#1\ldimen@i\relax
 \fi}
\def\dl@#1{\ifdl@ii\else
 \global\dl@iitrue \ldimen@ii\hunit \global\ldimen@ii#1\ldimen@ii\relax
 \fi}
\def\s@{\ifunshort@\else\global\short@true\fi}
\def\uns@{\ifshort@\else\global\unshort@true\global\short@false\fi}
\def\optioncodes@{\let\0\zero@\let\1\one@\let\a\a@\let\ds\ds@\let\dtX\dtX@
 \let\dtY\dtY@\let\da\da@\let\dx\dx@\let\dX\dX@\let\dY\dY@\let\dy\dy@
 \let\p\p@@\let\L\L@\let\l\l@\let\dL\dL@\let\dl\dl@\let\s\s@\let\uns\uns@}
\def\slopes@{\\161\\152\\143\\134\\255\\126\\357\\238\\349\\45{10}\\56{11}%
 \\11{12}\\65{13}\\54{14}\\43{15}\\32{16}\\53{17}\\21{18}\\52{19}\\31{20}%
 \\41{21}\\51{22}\\61{23}}
\newcount\tan@i
\newcount\tan@ip
\newcount\tan@ii
\newcount\tan@iip
\newdimen\slope@i
\newdimen\slope@ip
\newdimen\slope@ii
\newdimen\slope@iip
\newcount\angcount@
\newcount\extracount@
\def\slope@{{\slope@i=\secondy@ \advance\slope@i-\firsty@
 \ifN@\else\multiply\slope@i\m@ne\fi
 \slope@ii=\secondx@ \advance\slope@ii-\firstx@
 \ifE@\else\multiply\slope@ii\m@ne\fi
 \ifdim\slope@ii<\z@
  \global\tan@i6 \global\tan@ii\@ne \global\angcount@23
 \else
  \dimen@\slope@i \multiply\dimen@6
  \ifdim\dimen@<\slope@ii
   \global\tan@i\@ne \global\tan@ii6 \global\angcount@\@ne
  \else
   \dimen@\slope@ii \multiply\dimen@6
   \ifdim\dimen@<\slope@i
    \global\tan@i6 \global\tan@ii\@ne \global\angcount@23
   \else
    \tan@ip\z@ \tan@iip \@ne
    \def\\##1##2##3{\global\angcount@=##3\relax
     \slope@ip\slope@i \slope@iip\slope@ii
     \multiply\slope@iip##1\relax \multiply\slope@ip##2\relax
     \ifdim\slope@iip<\slope@ip
      \global\tan@ip=##1\relax \global\tan@iip=##2\relax
     \else
      \global\tan@i=##1\relax \global\tan@ii=##2\relax
      \def\\####1####2####3{}%
     \fi}%
    \slopes@
    \slope@i=\secondy@ \advance\slope@i-\firsty@
    \ifN@\else\multiply\slope@i\m@ne\fi
    \multiply\slope@i\tan@ii \multiply\slope@i\tan@iip \multiply\slope@i\tw@
    \count@\tan@i \multiply\count@\tan@iip
    \extracount@\tan@ip \multiply\extracount@\tan@ii
    \advance\count@\extracount@
    \slope@ii=\secondx@ \advance\slope@ii-\firstx@
    \ifE@\else\multiply\slope@ii\m@ne\fi
    \multiply\slope@ii\count@
    \ifdim\slope@i<\slope@ii
     \global\tan@i=\tan@ip \global\tan@ii=\tan@iip
     \global\advance\angcount@\m@ne
    \fi
   \fi
  \fi
 \fi}%
}
\def\slope@a#1{{\def\\##1##2##3{\ifnum##3=#1\global\tan@i=##1\relax
 \global\tan@ii=##2\relax\fi}\slopes@}}
\newcount\i@
\newcount\j@
\newcount\colcount@
\newcount\Colcount@
\newcount\tcolcount@
\newdimen\rowht@
\newdimen\rowdp@
\newcount\rowcount@
\newcount\Rowcount@
\newcount\maxcolrow@
\newtoks\colwidthtoks@
\newtoks\Rowheighttoks@
\newtoks\Rowdepthtoks@
\newtoks\widthtoks@
\newtoks\Widthtoks@
\newtoks\heighttoks@
\newtoks\Heighttoks@
\newtoks\depthtoks@
\newtoks\Depthtoks@
\newif\iffirstnewCDcr@
\def\dotoks@i{%
 \global\widthtoks@=\expandafter{\the\widthtoks@\else\getdim@\z@\fi}%
 \global\heighttoks@=\expandafter{\the\heighttoks@\else\getdim@\z@\fi}%
 \global\depthtoks@=\expandafter{\the\depthtoks@\else\getdim@\z@\fi}}
\def\dotoks@ii{%
 \global\widthtoks@{\ifcase\j@}%
 \global\heighttoks@{\ifcase\j@}%
 \global\depthtoks@{\ifcase\j@}}
\def\prenewCD@#1\endnewCD{\setbox\zer@
 \vbox{%
  \def\arrow@##1##2{{}}%
  \rowcount@\m@ne \colcount@\z@ \Colcount@\z@
  \firstnewCDcr@true \toks@{}%
  \widthtoks@{\ifcase\j@}%
  \Widthtoks@{\ifcase\i@}%
  \heighttoks@{\ifcase\j@}%
  \Heighttoks@{\ifcase\i@}%
  \depthtoks@{\ifcase\j@}%
  \Depthtoks@{\ifcase\i@}%
  \Rowheighttoks@{\ifcase\i@}%
  \Rowdepthtoks@{\ifcase\i@}%
  \Let@
  \everycr{%
   \noalign{%
    \global\advance\rowcount@\@ne
    \ifnum\colcount@<\Colcount@
    \else
     \global\Colcount@=\colcount@ \global\maxcolrow@=\rowcount@
    \fi
    \global\colcount@\z@
    \iffirstnewCDcr@
     \global\firstnewCDcr@false
    \else
     \edef\next@{\the\Rowheighttoks@\noexpand\or\noexpand\getdim@\the\rowht@}%
      \global\Rowheighttoks@=\expandafter{\next@}%
     \edef\next@{\the\Rowdepthtoks@\noexpand\or\noexpand\getdim@\the\rowdp@}%
      \global\Rowdepthtoks@=\expandafter{\next@}%
     \global\rowht@\z@ \global\rowdp@\z@
     \dotoks@i
     \edef\next@{\the\Widthtoks@\noexpand\or\the\widthtoks@}%
      \global\Widthtoks@=\expandafter{\next@}%
     \edef\next@{\the\Heighttoks@\noexpand\or\the\heighttoks@}%
      \global\Heighttoks@=\expandafter{\next@}%
     \edef\next@{\the\Depthtoks@\noexpand\or\the\depthtoks@}%
      \global\Depthtoks@=\expandafter{\next@}%
     \dotoks@ii
    \fi}}%
  \tabskip\z@
  \halign{&\setbox\zer@\hbox{\vrule height10\p@ width\z@ depth\z@
   $\m@th\displaystyle{##}$}\copy\zer@
   \ifdim\ht\zer@>\rowht@ \global\rowht@\ht\zer@ \fi
   \ifdim\dp\zer@>\rowdp@ \global\rowdp@\dp\zer@ \fi
   \global\advance\colcount@\@ne
   \edef\next@{\the\widthtoks@\noexpand\or\noexpand\getdim@\the\wd\zer@}%
    \global\widthtoks@=\expandafter{\next@}%
   \edef\next@{\the\heighttoks@\noexpand\or\noexpand\getdim@\the\ht\zer@}%
    \global\heighttoks@=\expandafter{\next@}%
   \edef\next@{\the\depthtoks@\noexpand\or\noexpand\getdim@\the\dp\zer@}%
    \global\depthtoks@=\expandafter{\next@}%
   \cr#1\crcr}}%
 \Rowcount@=\rowcount@
 \global\Widthtoks@=\expandafter{\the\Widthtoks@\fi\relax}%
 \edef\Width@##1##2{\i@=##1\relax\j@=##2\relax\the\Widthtoks@}%
 \global\Heighttoks@=\expandafter{\the\Heighttoks@\fi\relax}%
 \edef\Height@##1##2{\i@=##1\relax\j@=##2\relax\the\Heighttoks@}%
 \global\Depthtoks@=\expandafter{\the\Depthtoks@\fi\relax}%
 \edef\Depth@##1##2{\i@=##1\relax\j@=##2\relax\the\Depthtoks@}%
 \edef\next@{\the\Rowheighttoks@\noexpand\fi\relax}%
 \global\Rowheighttoks@=\expandafter{\next@}%
 \edef\Rowheight@##1{\i@=##1\relax\the\Rowheighttoks@}%
 \edef\next@{\the\Rowdepthtoks@\noexpand\fi\relax}%
 \global\Rowdepthtoks@=\expandafter{\next@}%
 \edef\Rowdepth@##1{\i@=##1\relax\the\Rowdepthtoks@}%
 \colwidthtoks@{\fi}%
 \setbox\zer@\vbox{%
  \unvbox\zer@
  \count@\rowcount@
  \loop
   \unskip\unpenalty
   \setbox\zer@\lastbox
   \ifnum\count@>\maxcolrow@ \advance\count@\m@ne
   \repeat
  \hbox{%
   \unhbox\zer@
   \count@\z@
   \loop
    \unskip
    \setbox\zer@\lastbox
    \edef\next@{\noexpand\or\noexpand\getdim@\the\wd\zer@\the\colwidthtoks@}%
     \global\colwidthtoks@=\expandafter{\next@}%
    \advance\count@\@ne
    \ifnum\count@<\Colcount@
    \repeat}}%
 \edef\next@{\noexpand\ifcase\noexpand\i@\the\colwidthtoks@}%
  \global\colwidthtoks@=\expandafter{\next@}%
 \edef\Colwidth@##1{\i@=##1\relax\the\colwidthtoks@}%
 \colwidthtoks@{}\Rowheighttoks@{}\Rowdepthtoks@{}\widthtoks@{}%
 \Widthtoks@{}\heighttoks@{}\Heighttoks@{}\depthtoks@{}\Depthtoks@{}%
}
\newcount\xoff@
\newcount\yoff@
\newcount\endcount@
\newcount\rcount@
\newdimen\firstx@
\newdimen\firsty@
\newdimen\secondx@
\newdimen\secondy@
\newdimen\tocenter@
\newdimen\charht@
\newdimen\charwd@
\def\outside@{\Err@{This arrow points outside the \string\newCD}}
\newif\ifsvertex@
\newif\iftvertex@
\def\arrow@#1#2{\xoff@=#1\relax\yoff@=#2\relax
 \count@\rowcount@ \advance\count@-\yoff@
 \ifnum\count@<\@ne \outside@ \else \ifnum\count@>\Rowcount@ \outside@ \fi\fi
 \count@\colcount@ \advance\count@\xoff@
 \ifnum\count@<\@ne \outside@ \else \ifnum\count@>\Colcount@ \outside@\fi\fi
 \tcolcount@\colcount@ \advance\tcolcount@\xoff@
 \Width@\rowcount@\colcount@ \tocenter@=-\getdim@ \divide\tocenter@\tw@
 \ifdim\getdim@=\z@
  \firstx@\z@ \firsty@\mathaxis@ \svertex@true
 \else
  \svertex@false
  \ifHshort@
   \Colwidth@\colcount@
    \ifE@ \firstx@=.5\getdim@ \else \firstx@=-.5\getdim@ \fi
  \else
   \ifE@ \firstx@=\getdim@ \else \firstx@=-\getdim@ \fi
   \divide\firstx@\tw@
  \fi
  \ifE@
   \ifH@ \advance\firstx@\thr@@\p@ \else \advance\firstx@-\thr@@\p@ \fi
  \else
   \ifH@ \advance\firstx@-\thr@@\p@ \else \advance\firstx@\thr@@\p@ \fi
  \fi
  \ifN@
   \Height@\rowcount@\colcount@ \firsty@=\getdim@
   \ifV@ \advance\firsty@\thr@@\p@ \fi
  \else
   \ifV@
    \Depth@\rowcount@\colcount@ \firsty@=-\getdim@
    \advance\firsty@-\thr@@\p@
   \else
    \firsty@\z@
   \fi
  \fi
 \fi
 \ifV@
 \else
  \Colwidth@\colcount@
  \ifE@ \secondx@=\getdim@ \else \secondx@=-\getdim@ \fi
  \divide\secondx@\tw@
  \ifE@ \else \getcgap@\colcount@ \advance\secondx@-\getdim@ \fi
  \endcount@=\colcount@ \advance\endcount@\xoff@
  \count@=\colcount@
  \ifE@
   \advance\count@\@ne
   \loop
    \ifnum\count@<\endcount@
    \Colwidth@\count@ \advance\secondx@\getdim@
    \getcgap@\count@ \advance\secondx@\getdim@
    \advance\count@\@ne
    \repeat
  \else
   \advance\count@\m@ne
   \loop
    \ifnum\count@>\endcount@
    \Colwidth@\count@ \advance\secondx@-\getdim@
    \getcgap@\count@ \advance\secondx@-\getdim@
    \advance\count@\m@ne
    \repeat
  \fi
  \Colwidth@\count@ \divide\getdim@\tw@
  \ifHshort@
  \else
   \ifE@ \advance\secondx@\getdim@ \else \advance\secondx@-\getdim@ \fi
  \fi
  \ifE@ \getcgap@\count@ \advance\secondx@\getdim@ \fi
  \rcount@\rowcount@ \advance\rcount@-\yoff@
  \Width@\rcount@\count@ \divide\getdim@\tw@
  \tvertex@false
  \ifH@\ifdim\getdim@=\z@\tvertex@true\Hshort@false\fi\fi
  \ifHshort@
  \else
   \ifE@ \advance\secondx@-\getdim@ \else \advance\secondx@\getdim@ \fi
  \fi
  \iftvertex@
   \advance\secondx@.4\p@
  \else
   \ifE@ \advance\secondx@-\thr@@\p@ \else \advance\secondx@\thr@@\p@ \fi
  \fi
 \fi
 \ifH@
 \else
  \ifN@
   \Rowheight@\rowcount@ \secondy@\getdim@
  \else
   \Rowdepth@\rowcount@ \secondy@-\getdim@
   \getrgap@\rowcount@ \advance\secondy@-\getdim@
  \fi
  \endcount@=\rowcount@ \advance\endcount@-\yoff@
  \count@=\rowcount@
  \ifN@
   \advance\count@\m@ne
   \loop
    \ifnum\count@>\endcount@
    \Rowheight@\count@ \advance\secondy@\getdim@
    \Rowdepth@\count@ \advance\secondy@\getdim@
    \getrgap@\count@ \advance\secondy@\getdim@
    \advance\count@\m@ne
    \repeat
  \else
   \advance\count@\@ne
   \loop
    \ifnum\count@<\endcount@
    \Rowheight@\count@ \advance\secondy@-\getdim@
    \Rowdepth@\count@ \advance\secondy@-\getdim@
    \getrgap@\count@ \advance\secondy@-\getdim@
    \advance\count@\@ne
    \repeat
  \fi
  \tvertex@false
  \ifV@\Width@\count@\colcount@\ifdim\getdim@=\z@\tvertex@true\fi\fi
  \ifN@
   \getrgap@\count@ \advance\secondy@\getdim@
   \Rowdepth@\count@ \advance\secondy@\getdim@
   \iftvertex@
    \advance\secondy@\mathaxis@
   \else
    \Depth@\count@\tcolcount@ \advance\secondy@-\getdim@
    \advance\secondy@-\thr@@\p@
   \fi
  \else
   \Rowheight@\count@ \advance\secondy@-\getdim@
   \iftvertex@
    \advance\secondy@\mathaxis@
   \else
    \Height@\count@\tcolcount@ \advance\secondy@\getdim@
    \advance\secondy@\thr@@\p@
   \fi
  \fi
 \fi
 \ifV@\else\advance\firstx@\sxdimen@\fi
 \ifH@\else\advance\firsty@\sydimen@\fi
 \iftX@
  \advance\secondy@\tXdimen@ii
  \advance\secondx@\tXdimen@i
  \slope@
 \else
  \iftY@
   \advance\secondy@\tYdimen@ii
   \advance\secondx@\tYdimen@i
   \slope@
   \secondy@=\secondx@ \advance\secondy@-\firstx@
   \ifNESW@ \else \multiply\secondy@\m@ne \fi
   \multiply\secondy@\tan@i \divide\secondy@\tan@ii \advance\secondy@\firsty@
  \else
   \ifa@
    \slope@
    \ifNESW@ \global\advance\angcount@\exacount@ \else
      \global\advance\angcount@-\exacount@ \fi
    \ifnum\angcount@>23 \angcount@23 \fi
    \ifnum\angcount@<\@ne \angcount@\@ne \fi
    \slope@a\angcount@
    \ifY@
     \advance\secondy@\Ydimen@
    \else
     \ifX@
      \advance\secondx@\Xdimen@
      \dimen@\secondx@ \advance\dimen@-\firstx@
      \ifNESW@\else\multiply\dimen@\m@ne\fi
      \multiply\dimen@\tan@i \divide\dimen@\tan@ii
      \advance\dimen@\firsty@ \secondy@=\dimen@
     \fi
    \fi
   \else
    \ifH@\else\ifV@\else\slope@\fi\fi
   \fi
  \fi
 \fi
 \ifH@\else\ifV@\else\ifsvertex@\else
  \dimen@=6\p@ \multiply\dimen@\tan@ii
  \count@=\tan@i \advance\count@\tan@ii \divide\dimen@\count@
  \ifE@ \advance\firstx@\dimen@ \else \advance\firstx@-\dimen@ \fi
  \multiply\dimen@\tan@i \divide\dimen@\tan@ii
  \ifN@ \advance\firsty@\dimen@ \else \advance\firsty@-\dimen@ \fi
 \fi\fi\fi
 \ifp@
  \ifH@\else\ifV@\else
   \getcos@\pdimen@ \advance\firsty@\dimen@ \advance\secondy@\dimen@
   \ifNESW@ \advance\firstx@-\dimen@ii \else \advance\firstx@\dimen@ii \fi
  \fi\fi
 \fi
 \ifH@\else\ifV@\else
  \ifnum\tan@i>\tan@ii
   \charht@=10\p@ \charwd@=10\p@
   \multiply\charwd@\tan@ii \divide\charwd@\tan@i
  \else
   \charwd@=10\p@ \charht@=10\p@
   \divide\charht@\tan@ii \multiply\charht@\tan@i
  \fi
  \ifnum\tcount@=\thr@@
   \ifN@ \advance\secondy@-.3\charht@ \else\advance\secondy@.3\charht@ \fi
  \fi
  \ifnum\scount@=\tw@
   \ifE@ \advance\firstx@.3\charht@ \else \advance\firstx@-.3\charht@ \fi
  \fi
  \ifnum\tcount@=12
   \ifN@ \advance\secondy@-\charht@ \else \advance\secondy@\charht@ \fi
  \fi
  \iftY@
  \else
   \ifa@
    \ifX@
    \else
     \secondx@\secondy@ \advance\secondx@-\firsty@
     \ifNESW@\else\multiply\secondx@\m@ne\fi
     \multiply\secondx@\tan@ii \divide\secondx@\tan@i
     \advance\secondx@\firstx@
    \fi
   \fi
  \fi
 \fi\fi
 \ifH@\harrow@\else\ifV@\varrow@\else\arrow@@\fi\fi}
\newdimen\mathaxis@
\mathaxis@90\p@ \divide\mathaxis@36
\def\harrow@b{\ifE@\hskip\tocenter@\hskip\firstx@\fi}
\def\harrow@bb{\ifE@\hskip\xdimen@\else\hskip\Xdimen@\fi}
\def\harrow@e{\ifE@\else\hskip-\firstx@\hskip-\tocenter@\fi}
\def\harrow@ee{\ifE@\hskip-\Xdimen@\else\hskip-\xdimen@\fi}
\def\harrow@{\dimen@\secondx@\advance\dimen@-\firstx@
 \ifE@ \let\next@\rlap \else  \multiply\dimen@\m@ne \let\next@\llap \fi
 \next@{%
  \harrow@b
  \smash{\raise\pdimen@\hbox to\dimen@
   {\harrow@bb\arrow@ii
    \ifnum\arrcount@=\m@ne \else \ifnum\arrcount@=\thr@@ \else
     \ifE@
      \ifnum\scount@=\m@ne
      \else
       \ifcase\scount@\or\or\char118 \or\char117 \or\or\or\char119 \or
       \char120 \or\char121 \or\char122 \or\or\or\arrow@i\char125 \or
       \char117 \hskip\thr@@\p@\char117 \hskip-\thr@@\p@\fi
      \fi
     \else
      \ifnum\tcount@=\m@ne
      \else
       \ifcase\tcount@\char117 \or\or\char117 \or\char118 \or\char119 \or
       \char120\or\or\or\or\or\char121 \or\char122 \or\arrow@i\char125
       \or\char117 \hskip\thr@@\p@\char117 \hskip-\thr@@\p@\fi
      \fi
     \fi
    \fi\fi
    \dimen@\mathaxis@ \advance\dimen@.2\p@
    \dimen@ii\mathaxis@ \advance\dimen@ii-.2\p@
    \ifnum\arrcount@=\m@ne
     \let\leads@\null
    \else
     \ifcase\arrcount@
      \def\leads@{\hrule height\dimen@ depth-\dimen@ii}\or
      \def\leads@{\hrule height\dimen@ depth-\dimen@ii}\or
      \def\leads@{\hbox to10\p@{%
       \leaders\hrule height\dimen@ depth-\dimen@ii\hfil
       \hfil
      \leaders\hrule height\dimen@ depth-\dimen@ii\hskip\z@ plus2fil\relax
       \hfil
       \leaders\hrule height\dimen@ depth-\dimen@ii\hfil}}\or
     \def\leads@{\hbox{\hbox to10\p@{\dimen@\mathaxis@ \advance\dimen@1.2\p@
       \dimen@ii\dimen@ \advance\dimen@ii-.4\p@
       \leaders\hrule height\dimen@ depth-\dimen@ii\hfil}%
       \kern-10\p@
       \hbox to10\p@{\dimen@\mathaxis@ \advance\dimen@-1.2\p@
       \dimen@ii\dimen@ \advance\dimen@ii-.4\p@
       \leaders\hrule height\dimen@ depth-\dimen@ii\hfil}}}\fi
    \fi
    \cleaders\leads@\hfil
    \ifnum\arrcount@=\m@ne\else\ifnum\arrcount@=\thr@@\else
     \arrow@i
     \ifE@
      \ifnum\tcount@=\m@ne
      \else
       \ifcase\tcount@\char119 \or\or\char119 \or\char120 \or\char121 \or
       \char122 \or \or\or\or\or\char123\or\char124 \or
       \char125 \or\char119 \hskip-\thr@@\p@\char119 \hskip\thr@@\p@\fi
      \fi
     \else
      \ifcase\scount@\or\or\char120 \or\char119 \or\or\or\char121 \or\char122
      \or\char123 \or\char124 \or\or\or\char125 \or
      \char119 \hskip-\thr@@\p@\char119 \hskip\thr@@\p@\fi
     \fi
    \fi\fi
    \harrow@ee}}%
  \harrow@e}%
 \iflabel@i
  \dimen@ii\z@ \setbox\zer@\hbox{$\m@th\tsize@@\label@i$}%
  \ifnum\arrcount@=\m@ne
  \else
   \advance\dimen@ii\mathaxis@
   \advance\dimen@ii\dp\zer@ \advance\dimen@ii\tw@\p@
   \ifnum\arrcount@=\thr@@ \advance\dimen@ii\tw@\p@ \fi
  \fi
  \advance\dimen@ii\pdimen@
  \next@{\harrow@b\smash{\raise\dimen@ii\hbox to\dimen@
   {\harrow@bb\hskip\tw@\ldimen@i\hfil\box\zer@\hfil\harrow@ee}}\harrow@e}%
 \fi
 \iflabel@ii
  \ifnum\arrcount@=\m@ne
  \else
   \setbox\zer@\hbox{$\m@th\tsize@\label@ii$}%
   \dimen@ii-\ht\zer@ \advance\dimen@ii-\tw@\p@
   \ifnum\arrcount@=\thr@@ \advance\dimen@ii-\tw@\p@ \fi
   \advance\dimen@ii\mathaxis@ \advance\dimen@ii\pdimen@
   \next@{\harrow@b\smash{\raise\dimen@ii\hbox to\dimen@
    {\harrow@bb\hskip\tw@\ldimen@ii\hfil\box\zer@\hfil\harrow@ee}}\harrow@e}%
  \fi
 \fi}
\let\tsize@\tsize
\def\tsizenewCDlabels{\let\tsize@\tsize}
\def\ssizenewCDlabels{\let\tsize@\ssize}
\def\tsize@@{\ifnum\arrcount@=\m@ne\else\tsize@\fi}
\def\varrow@{\dimen@\secondy@ \advance\dimen@-\firsty@
 \ifN@ \else \multiply\dimen@\m@ne \fi
 \setbox\zer@\vbox to\dimen@
  {\ifN@ \vskip-\Ydimen@ \else \vskip\ydimen@ \fi
   \ifnum\arrcount@=\m@ne\else\ifnum\arrcount@=\thr@@\else
    \hbox{\arrow@iii
     \ifN@
      \ifnum\tcount@=\m@ne
      \else
       \ifcase\tcount@\char117 \or\or\char117 \or\char118 \or\char119 \or
       \char120 \or\or\or\or\or\char121 \or\char122 \or\char123 \or
       \vbox{\hbox{\char117 }\nointerlineskip\vskip\thr@@\p@
       \hbox{\char117 }\vskip-\thr@@\p@}\fi
      \fi
     \else
      \ifcase\scount@\or\or\char118 \or\char117 \or\or\or\char119 \or
      \char120 \or\char121 \or\char122 \or\or\or\char123 \or
      \vbox{\hbox{\char117 }\nointerlineskip\vskip\thr@@\p@
      \hbox{\char117 }\vskip-\thr@@\p@}\fi
     \fi}%
    \nointerlineskip
   \fi\fi
   \ifnum\arrcount@=\m@ne
    \let\leads@\null
   \else
    \ifcase\arrcount@\let\leads@\vrule\or\let\leads@\vrule\or
    \def\leads@{\vbox to10\p@{%
     \hrule height 1.67\p@ depth\z@ width.4\p@
     \vfil
     \hrule height 3.33\p@ depth\z@ width.4\p@
     \vfil
     \hrule height 1.67\p@ depth\z@ width.4\p@}}\or
    \def\leads@{\hbox{\vrule height\p@\hskip\tw@\p@\vrule}}\fi
   \fi
  \cleaders\leads@\vfill\nointerlineskip
   \ifnum\arrcount@=\m@ne\else\ifnum\arrcount@=\thr@@\else
    \hbox{\arrow@iv
     \ifN@
      \ifcase\scount@\or\or\char118 \or\char117 \or\or\or\char119 \or
      \char120 \or\char121 \or\char122 \or\or\or\arrow@iii\char123 \or
      \vbox{\hbox{\char117 }\nointerlineskip\vskip-\thr@@\p@
      \hbox{\char117 }\vskip\thr@@\p@}\fi
     \else
      \ifnum\tcount@=\m@ne
      \else
       \ifcase\tcount@\char117 \or\or\char117 \or\char118 \or\char119 \or
       \char120 \or\or\or\or\or\char121 \or\char122 \or\arrow@iii\char123 \or
       \vbox{\hbox{\char117 }\nointerlineskip\vskip-\thr@@\p@
       \hbox{\char117 }\vskip\thr@@\p@}\fi
      \fi
     \fi}%
   \fi\fi
   \ifN@\vskip\ydimen@\else\vskip-\Ydimen@\fi}%
 \ifN@
  \dimen@ii\firsty@
 \else
  \dimen@ii-\firsty@ \advance\dimen@ii\ht\zer@ \multiply\dimen@ii\m@ne
 \fi
 \rlap{\smash{\hskip\tocenter@ \hskip\pdimen@ \raise\dimen@ii \box\zer@}}%
 \iflabel@i
  \setbox\zer@\vbox to\dimen@{\vfil
   \hbox{$\m@th\tsize@@\label@i$}\vskip\tw@\ldimen@i\vfil}%
  \rlap{\smash{\hskip\tocenter@ \hskip\pdimen@
  \ifnum\arrcount@=\m@ne \let\next@\relax \else \let\next@\llap \fi
  \next@{\raise\dimen@ii\hbox{\ifnum\arrcount@=\m@ne \hskip-.5\wd\zer@ \fi
   \box\zer@ \ifnum\arrcount@=\m@ne \else \hskip\tw@\p@ \fi}}}}%
 \fi
 \iflabel@ii
  \ifnum\arrcount@=\m@ne
  \else
   \setbox\zer@\vbox to\dimen@{\vfil
    \hbox{$\m@th\tsize@\label@ii$}\vskip\tw@\ldimen@ii\vfil}%
   \rlap{\smash{\hskip\tocenter@ \hskip\pdimen@
   \rlap{\raise\dimen@ii\hbox{\ifnum\arrcount@=\thr@@ \hskip4.5\p@ \else
    \hskip2.5\p@ \fi\box\zer@}}}}%
  \fi
 \fi
}
\newdimen\goal@
\newdimen\shifted@
\newcount\Tcount@
\newcount\Scount@
\newbox\shaft@
\newcount\slcount@
\def\getcos@#1{%
 \ifnum\tan@i<\tan@ii
  \dimen@#1%
  \ifnum\slcount@<8 \count@9 \else \ifnum\slcount@<12 \count@8 \else
   \count@7 \fi\fi
  \multiply\dimen@\count@ \divide\dimen@10
  \dimen@ii\dimen@ \multiply\dimen@ii\tan@i \divide\dimen@ii\tan@ii
 \else
  \dimen@ii#1%
  \count@-\slcount@ \advance\count@24
  \ifnum\count@<8 \count@9 \else \ifnum\count@<12 \count@8
   \else\count@7 \fi\fi
  \multiply\dimen@ii\count@ \divide\dimen@ii10
  \dimen@\dimen@ii \multiply\dimen@\tan@ii \divide\dimen@\tan@i
 \fi}
\newdimen\adjust@
\def\Nnext@{\ifN@\let\next@\raise\else\let\next@\lower\fi}
\def\arrow@@{\slcount@\angcount@
 \ifNESW@
  \ifnum\angcount@<10
   \let\arrowfont@=\arrow@i \advance\angcount@\m@ne \multiply\angcount@13
  \else
   \ifnum\angcount@<19
    \let\arrowfont@=\arrow@ii \advance\angcount@-10 \multiply\angcount@13
   \else
    \let\arrowfont@=\arrow@iii \advance\angcount@-19 \multiply\angcount@13
  \fi\fi
  \Tcount@\angcount@
 \else
  \ifnum\angcount@<5
   \let\arrowfont@=\arrow@iii \advance\angcount@\m@ne \multiply\angcount@13
   \advance\angcount@65
  \else
   \ifnum\angcount@<14
    \let\arrowfont@=\arrow@iv \advance\angcount@-5 \multiply\angcount@13
   \else
    \ifnum\angcount@<23
     \let\arrowfont@=\arrow@v \advance\angcount@-14 \multiply\angcount@13
    \else
     \let\arrowfont@=\arrow@i \angcount@=117
  \fi\fi\fi
  \ifnum\angcount@=117 \Tcount@=115 \else\Tcount@\angcount@ \fi
 \fi
 \Scount@\Tcount@
 \ifE@
  \ifnum\tcount@=\z@ \advance\Tcount@\tw@ \else\ifnum\tcount@=13
   \advance\Tcount@\tw@ \else \advance\Tcount@\tcount@ \fi\fi
  \ifnum\scount@=\z@ \else \ifnum\scount@=13 \advance\Scount@\thr@@ \else
   \advance\Scount@\scount@ \fi\fi
 \else
  \ifcase\tcount@\advance\Tcount@\thr@@\or\or\advance\Tcount@\thr@@\or
  \advance\Tcount@\tw@\or\advance\Tcount@6 \or\advance\Tcount@7
  \or\or\or\or\or \advance\Tcount@8 \or\advance\Tcount@9 \or
  \advance\Tcount@12 \or\advance\Tcount@\thr@@\fi
  \ifcase\scount@\or\or\advance\Scount@\thr@@\or\advance\Scount@\tw@\or
  \or\or\advance\Scount@4 \or\advance\Scount@5 \or\advance\Scount@10
  \or\advance\Scount@11 \or\or\or\advance\Scount@12 \or\advance
  \Scount@\tw@\fi
 \fi
 \ifcase\arrcount@\or\or\advance\angcount@\@ne\else\fi
 \ifN@ \shifted@=\firsty@ \else\shifted@=-\firsty@ \fi
 \ifE@ \else\advance\shifted@\charht@ \fi
 \goal@=\secondy@ \advance\goal@-\firsty@
 \ifN@\else\multiply\goal@\m@ne\fi
 \setbox\shaft@\hbox{\arrowfont@\char\angcount@}%
 \ifnum\arrcount@=\thr@@
  \getcos@{1.5\p@}%
  \setbox\shaft@\hbox to\wd\shaft@{\arrowfont@
   \rlap{\hskip\dimen@ii
    \smash{\ifNESW@\let\next@\lower\else\let\next@\raise\fi
     \next@\dimen@\hbox{\arrowfont@\char\angcount@}}}%
   \rlap{\hskip-\dimen@ii
    \smash{\ifNESW@\let\next@\raise\else\let\next@\lower\fi
      \next@\dimen@\hbox{\arrowfont@\char\angcount@}}}\hfil}%
 \fi
 \rlap{\smash{\hskip\tocenter@\hskip\firstx@
  \ifnum\arrcount@=\m@ne
  \else
   \ifnum\arrcount@=\thr@@
   \else
    \ifnum\scount@=\m@ne
    \else
     \ifnum\scount@=\z@
     \else
      \setbox\zer@\hbox{\ifnum\angcount@=117 \arrow@v\else\arrowfont@\fi
       \char\Scount@}%
      \ifNESW@
       \ifnum\scount@=\tw@
        \dimen@=\shifted@ \advance\dimen@-\charht@
        \ifN@\hskip-\wd\zer@\fi
        \Nnext@
        \next@\dimen@\copy\zer@
        \ifN@\else\hskip-\wd\zer@\fi
       \else
        \Nnext@
        \ifN@\else\hskip-\wd\zer@\fi
        \next@\shifted@\copy\zer@
        \ifN@\hskip-\wd\zer@\fi
       \fi
       \ifnum\scount@=12
        \advance\shifted@\charht@ \advance\goal@-\charht@
        \ifN@ \hskip\wd\zer@ \else \hskip-\wd\zer@ \fi
       \fi
       \ifnum\scount@=13
        \getcos@{\thr@@\p@}%
        \ifN@ \hskip\dimen@ \else \hskip-\wd\zer@ \hskip-\dimen@ \fi
        \adjust@\shifted@ \advance\adjust@\dimen@ii
        \Nnext@
        \next@\adjust@\copy\zer@
        \ifN@ \hskip-\dimen@ \hskip-\wd\zer@ \else \hskip\dimen@ \fi
       \fi
      \else
       \ifN@\hskip-\wd\zer@\fi
       \ifnum\scount@=\tw@
        \ifN@ \hskip\wd\zer@ \else \hskip-\wd\zer@ \fi
        \dimen@=\shifted@ \advance\dimen@-\charht@
        \Nnext@
        \next@\dimen@\copy\zer@
        \ifN@\hskip-\wd\zer@\fi
       \else
        \Nnext@
        \next@\shifted@\copy\zer@
        \ifN@\else\hskip-\wd\zer@\fi
       \fi
       \ifnum\scount@=12
        \advance\shifted@\charht@ \advance\goal@-\charht@
        \ifN@ \hskip-\wd\zer@ \else \hskip\wd\zer@ \fi
       \fi
       \ifnum\scount@=13
        \getcos@{\thr@@\p@}%
        \ifN@ \hskip-\wd\zer@ \hskip-\dimen@ \else \hskip\dimen@ \fi
        \adjust@\shifted@ \advance\adjust@\dimen@ii
        \Nnext@
        \next@\adjust@\copy\zer@
        \ifN@ \hskip\dimen@ \else \hskip-\dimen@ \hskip-\wd\zer@ \fi
       \fi	
      \fi
  \fi\fi\fi\fi
  \ifnum\arrcount@=\m@ne
  \else
   \loop
    \ifdim\goal@>\charht@
    \ifE@\else\hskip-\charwd@\fi
    \Nnext@
    \next@\shifted@\copy\shaft@
    \ifE@\else\hskip-\charwd@\fi
    \advance\shifted@\charht@ \advance\goal@ -\charht@
    \repeat
   \ifdim\goal@>\z@
    \dimen@=\charht@ \advance\dimen@-\goal@
    \divide\dimen@\tan@i \multiply\dimen@\tan@ii
    \ifE@ \hskip-\dimen@ \else \hskip-\charwd@ \hskip\dimen@ \fi
    \adjust@=\shifted@ \advance\adjust@-\charht@ \advance\adjust@\goal@
    \Nnext@
    \next@\adjust@\copy\shaft@
    \ifE@ \else \hskip-\charwd@ \fi
   \else
    \adjust@=\shifted@ \advance\adjust@-\charht@
   \fi
  \fi
  \ifnum\arrcount@=\m@ne
  \else
   \ifnum\arrcount@=\thr@@
   \else
    \ifnum\tcount@=\m@ne
    \else
     \setbox\zer@
      \hbox{\ifnum\angcount@=117 \arrow@v\else\arrowfont@\fi\char\Tcount@}%
     \ifnum\tcount@=\thr@@
      \advance\adjust@\charht@
      \ifE@\else\ifN@\hskip-\charwd@\else\hskip-\wd\zer@\fi\fi
     \else
      \ifnum\tcount@=12
       \advance\adjust@\charht@
       \ifE@\else\ifN@\hskip-\charwd@\else\hskip-\wd\zer@\fi\fi
      \else
       \ifE@\hskip-\wd\zer@\fi
     \fi\fi
     \Nnext@
     \next@\adjust@\copy\zer@
     \ifnum\tcount@=13
      \hskip-\wd\zer@
      \getcos@{\thr@@\p@}%
      \ifE@\hskip-\dimen@ \else\hskip\dimen@ \fi
      \advance\adjust@-\dimen@ii
      \Nnext@
      \next@\adjust@\box\zer@
     \fi
  \fi\fi\fi}}%
 \iflabel@i
  \rlap{\hskip\tocenter@
  \dimen@\firstx@ \advance\dimen@\secondx@ \divide\dimen@\tw@
  \advance\dimen@\ldimen@i
  \dimen@ii\firsty@ \advance\dimen@ii\secondy@ \divide\dimen@ii\tw@
  \multiply\ldimen@i\tan@i \divide\ldimen@i\tan@ii
  \ifNESW@ \advance\dimen@ii\ldimen@i \else \advance\dimen@ii-\ldimen@i \fi
  \setbox\zer@\hbox{\ifNESW@\else\ifnum\arrcount@=\thr@@\hskip4\p@\else
   \hskip\tw@\p@\fi\fi
   $\m@th\tsize@@\label@i$\ifNESW@\ifnum\arrcount@=\thr@@\hskip4\p@\else
   \hskip\tw@\p@\fi\fi}%
  \ifnum\arrcount@=\m@ne
   \ifNESW@ \advance\dimen@.5\wd\zer@ \advance\dimen@\p@ \else
    \advance\dimen@-.5\wd\zer@ \advance\dimen@-\p@ \fi
   \advance\dimen@ii-.5\ht\zer@
  \else
   \advance\dimen@ii\dp\zer@
   \ifnum\slcount@<6 \advance\dimen@ii\tw@\p@ \fi
  \fi
  \hskip\dimen@
  \ifNESW@ \let\next@\llap \else\let\next@\rlap \fi
  \next@{\smash{\raise\dimen@ii\box\zer@}}}%
 \fi
 \iflabel@ii
  \ifnum\arrcount@=\m@ne
  \else
   \rlap{\hskip\tocenter@
   \dimen@\firstx@ \advance\dimen@\secondx@ \divide\dimen@\tw@
   \ifNESW@ \advance\dimen@\ldimen@ii \else \advance\dimen@-\ldimen@ii \fi
   \dimen@ii\firsty@ \advance\dimen@ii\secondy@ \divide\dimen@ii\tw@
   \multiply\ldimen@ii\tan@i \divide\ldimen@ii\tan@ii
   \advance\dimen@ii\ldimen@ii
   \setbox\zer@\hbox{\ifNESW@\ifnum\arrcount@=\thr@@\hskip4\p@\else
    \hskip\tw@\p@\fi\fi
    $\m@th\tsize@\label@ii$\ifNESW@\else\ifnum\arrcount@=\thr@@\hskip4\p@
    \else\hskip\tw@\p@\fi\fi}%
   \advance\dimen@ii-\ht\zer@
   \ifnum\slcount@<9 \advance\dimen@ii-\thr@@\p@ \fi
   \ifNESW@ \let\next@\rlap \else \let\next@\llap \fi
   \hskip\dimen@\next@{\smash{\raise\dimen@ii\box\zer@}}}%
  \fi
 \fi
}
\def\outnewCD@#1{\def#1{\Err@{\string#1 must not be used within \string\newCD}}}
\newskip\prenewCDskip@
\newskip\postnewCDskip@
\prenewCDskip@\z@
\postnewCDskip@\z@
\def\prenewCDspace#1{\RIfMIfI@
 \onlydmatherr@\prenewCDspace\else\advance\prenewCDskip@#1\relax\fi\else
 \onlydmatherr@\prenewCDspace\fi}
\def\postnewCDspace#1{\RIfMIfI@
 \onlydmatherr@\postnewCDspace\else\advance\postnewCDskip@#1\relax\fi\else
 \onlydmatherr@\postnewCDspace\fi}
\def\predisplayspace#1{\RIfMIfI@
 \onlydmatherr@\predisplayspace\else
 \advance\abovedisplayskip#1\relax
 \advance\abovedisplayshortskip#1\relax\fi
 \else\onlydmatherr@\prenewCDspace\fi}
\def\postdisplayspace#1{\RIfMIfI@
 \onlydmatherr@\postdisplayspace\else
 \advance\belowdisplayskip#1\relax
 \advance\belowdisplayshortskip#1\relax\fi
 \else\onlydmatherr@\postdisplayspace\fi}
\def\PrenewCDSpace#1{\global\prenewCDskip@#1\relax}
\def\PostnewCDSpace#1{\global\postnewCDskip@#1\relax}
\def\newCD#1\endnewCD{%
 \outnewCD@\cgaps\outnewCD@\rgaps\outnewCD@\Cgaps\outnewCD@\Rgaps
 \prenewCD@#1\endnewCD
 \advance\abovedisplayskip\prenewCDskip@
 \advance\abovedisplayshortskip\prenewCDskip@
 \advance\belowdisplayskip\postnewCDskip@
 \advance\belowdisplayshortskip\postnewCDskip@
 \vcenter{\vskip\prenewCDskip@ \Let@ \colcount@\@ne \rowcount@\z@
  \everycr{%
   \noalign{%
    \ifnum\rowcount@=\Rowcount@
    \else
     \global\nointerlineskip
     \getrgap@\rowcount@ \vskip\getdim@
     \global\advance\rowcount@\@ne \global\colcount@\@ne
    \fi}}%
  \tabskip\z@
  \halign{&\global\xoff@\z@ \global\yoff@\z@
   \getcgap@\colcount@ \hskip\getdim@
   \hfil\vrule height10\p@ width\z@ depth\z@
   $\m@th\displaystyle{##}$\hfil
   \global\advance\colcount@\@ne\cr
   #1\crcr}\vskip\postnewCDskip@}%
 \prenewCDskip@\z@\postnewCDskip@\z@
 \def\getcgap@##1{\ifcase##1\or\getdim@\z@\else\getdim@\standardcgap\fi}%
 \def\getrgap@##1{\ifcase##1\getdim@\z@\else\getdim@\standardrgap\fi}%
 \let\Width@\relax\let\Height@\relax\let\Depth@\relax\let\Rowheight@\relax
 \let\Rowdepth@\relax\let\Colwdith@\relax
}
\catcode`\@=\active
%\endinput
%end of diag.tex
\hsize 30pc
\vsize 47pc
\def\nmb#1#2{#2}         % used for renumbering, TeX should ignore.
\def\cit#1#2{\ifx#1!\cite{#2}\else#2\fi} %for citing references
\def\totoc{}             %= to table of content, invoked by kms-book.sty
\def\idx{}               % for producing index, invoked by kms-book.sty
\def\ign#1{}             %=ignore, invisible entry for the index only

\redefine\o{\circ}

\define\al{\alpha}
\define\be{\beta}

\define\de{\delta}
\define\ep{\varepsilon}

\define\ka{\kappa}
\define\la{\lambda}

\define\ph{\varphi}

\define\om{\omega}

\define\La{\Lambda}

\redefine\i{^{-1}}
\define\row#1#2#3{#1_{#2},\ldots,#1_{#3}}
\define\x{\times}

\define\ev{\operatorname{ev}}
\define\Id{\operatorname{Id}}
\predefine\LL{\L}
\redefine\L{{\Cal L}}

\define\Mf{\Cal Mf}

\define\pr{\operatorname{pr}}
\define\Hom{\operatorname{Hom}}

\def\today{\ifcase\month\or
 January\or February\or March\or April\or May\or June\or
 July\or August\or September\or October\or November\or December\fi
 \space\number\day, \number\year}
\topmatter
\title Product preserving functors of infinite dimensional manifolds 
\endtitle
\author  Andreas Kriegl\\
Peter W. Michor  \endauthor
\affil
Institut f\"ur Mathematik, Universit\"at Wien,\\
Strudlhofgasse 4, A-1090 Wien, Austria.\\
Erwin Schr\"odinger Institut f\"ur Mathematische Physik,
Boltzmanngasse 9, A-1090 Wien, Austria
\endaffil
\address
Institut f\"ur Mathematik, Universit\"at Wien,
Strudlhofgasse 4, A-1090 Wien, Austria
\endaddress
\email kriegl\@pap.univie.ac.at,  Peter.Michor\@esi.ac.at \endemail
\dedicatory 
For Ivan Kol\'a\v r, on the occasion of his 60th birthday.
\enddedicatory
\date {\today} \enddate
\thanks 
Supported by `Fonds zur F\"orderung der wissenscahftlichen 
Forschung, Projekt P~10037~PHY'.
\endthanks
\abstract 
The theory of product preserving functors and Weil functors is 
partly extended to infinite dimensional manifolds, using the theory 
of $C^\infty$-algebras. 
\endabstract
\endtopmatter
%\input amspptb.sty
%\userunningheads
%\def\leftheadtext{\smc Andreas Kriegl, Peter W. Michor}
%\def\rightheadtext{ Product preserving functors}
%\def\bottremark{\today\hfill}

\document

\heading Table of contents \endheading
%\input \jobname.toc
%\loadtoc
%\loadindex
\noindent 1. Introduction \leaders \hbox to 
1em{\hss .\hss }\hfill {\eightrm 1}\par  
\noindent 2. Infinite dimensional manifolds \leaders \hbox to 
1em{\hss .\hss }\hfill {\eightrm 2}\par  
\noindent 3. Weil functors on infinite dimensional manifolds 
\leaders \hbox to 1em{\hss .\hss }\hfill {\eightrm 3}\par  
\noindent 4. Product preserving functors from finite dimensional 
manifolds \par to infinite dimensional ones \leaders \hbox to 
1em{\hss .\hss }\hfill {\eightrm 10}\par  

\head\totoc\nmb0{1}. Introduction \endhead

In competition to the theory of jets of Ehresmann Andr\'e Weil in 
\cit!{21} explained a construction which is on the one hand more 
restricted than that of jets since it allows only for covariant 
constructions in the sense of category theory (contravariant in the 
sense of differential geometry), but is more flexible since it uses 
more input: a finite dimensional formally real algebra. Later it was 
realized that Weil's construction describes all product preserving 
bundle functors on the category of finite dimension manifolds. This 
was developed independently by \cit!{1}, \cit!{4}, less completely 
by \cit!{13}, and exposed in detail in chapter 
VIII of \cit!{6}. A jet-like approach to Weil's construction was 
given in \cit!{18} and used by Kol\'a\v r in \cit!{5} to discuss natural 
transformations.
The purpose of this paper is to present some versions of this theory  
in the realm of infinite dimensional manifolds, in the setting of 
convenient calculus as in \cit!{2}, \cit!{10}, \cit!{8}. 

\head\totoc\nmb0{2}. Infinite dimensional manifolds \endhead

\subhead\nmb.{2.1}. Calculus in infinite dimensions \endsubhead
A locally convex vector space $E$ is called \idx{\it convenient} if 
for each smooth curve the Riemann integrals over compact intervals 
exist (this is a weak completeness condition). The final topology on 
$E$ with respect to all smooth curves is called the 
$c^\infty$-topology.
A mapping between ($c^\infty$-) open subsets of convenient vector spaces 
is called smooth if it maps smooth curves to smooth curves. 
Multilinear mappings are smooth if and only if they are bounded. This 
gives a meaningful theory which up to Fr\'echet spaces coincides with 
any reasonable theory of smooth mappings. 
The main additional property is cartesian closedness: If $U$, $V$ are 
$c^\infty$-open subsets in and if $G$ is a convenient vector space, 
then $C^\infty(V,G)$ is again a convenient vector space (with the 
locally convex topology of convergence of 
compositions with smooth curves in $V$, uniformly on compact 
intervals, in all derivatives separately), and we have  
$$
C^\infty(U,C^\infty(V,G))\cong C^\infty(U\x V, G).
$$
Expositions of this theory 
can be found in \cit!{2}, \cit!{10}, \cit!{8}, e.g.
Real analytic and holomorphic versions of this theory are also 
available, \cit!{11}, \cit!{7}, \cit!{10}.

\subhead \nmb.{2.2}. Manifolds \endsubhead 
A chart $(U,u)$ on a set $M$ is a bijection $u:U\to
u(U)\subseteq E_U$ from a subset $U\subseteq M$ onto a
$c^\infty$-open subset of a convenient vector space $E_U$. 

For two charts $(U_\al,u_\al)$ and $(U_\be,u_\be)$ on  $M$ the mapping 
$u_{\al\be} := u_\al\o u_\be\i:u_\be(U_{\al\be}) \to
u_\al(U_{\al\be})$ for $\al$, $\be\in A$
is called the \idx{\it chart changing\ign{ mapping}}, where 
$U_{\al\be} := U_\al\cap U_\be$.
A family $(U_\al,u_\al)_{\al\in A}$ of charts on
$M$ is called
an \idx{\it atlas} for $M$, if 
the $U_\al$ form a cover of $M$ and all chart changings
$u_{\al\be}$ are defined on $c^\infty$-open subsets. 

An atlas $(U_\al,u_\al)_{\al\in A}$ for $M$ is said to be a 
$C^\infty$-\idx{\it\ign{smooth }atlas}, if all chart changings 
$u_{\al\be}:u_\be(U_{\al\be}) \to u_\al(U_{\al\be})$ are smooth. Two 
$C^\infty$-atlases are called $C^\infty$-\idx{\it equivalent\ign{ 
atlases}}, if their union is again a $C^\infty$-atlas for $M$. An 
equivalence class of $C^\infty$-atlases is called a 
$C^\infty$\idx{\it-structure} on $M$. The union of all atlases in an 
equivalence class is again an atlas, the \idx{\it maximal} atlas for 
this $C^\infty$-structure. A $C^\infty$-\idx{\it manifold} $M$ is a 
set together with a $C^\infty$-structure on it.  

\subheading{\nmb.{2.3}} A mapping $f: M \to N$ between manifolds is
called smooth if for each $x\in M$ and each chart $(V,v)$ on 
$N$ with $f(x)\in V$ there is a chart $(U,u)$ on $M$ with $x\in U$,
$f(U)\subseteq V$, such that $v\o f\o u\i$ is smooth. This is the 
case if and only if $f\o c$ is smooth for each smooth curve 
$c:\Bbb R\to M$.

We will denote by $C^\infty(M,N)$ the space of all
$C^\infty$-mappings from $M$ to $N$. 

\subhead\nmb.{2.4}. The topology of a manifold \endsubhead
The \idx{\it natural topology\ign{ on a manifold}} on a manifold $M$ is the 
identification topology with respect to some (smooth) atlas
$(u_\al:M\supseteq U_\al\to u_\al(U_\al)\subseteq E_\al)$, 
where a subset
$W\subseteq M$ is open if and only if $u_\al(U_\al\cap W)$ is
$c^\infty$-open in $E_\al$ for all $\al$.
This topology depends only on the structure, since
diffeomorphisms are homeomorphisms for the 
$c^\infty$-topologies.
It is also the final topology with
respect to all inverses of chart mappings in one atlas. 
It is also the final topology with respect to all smooth curves.
For a (smooth) manifold we will require certain properties for
the natural topology, which will be specified when needed, like: 
\roster
\item \idx{\it Smoothly Hausdorff}: The smooth functions in 
       $C^\infty(M,\Bbb R)$ separate points in $M$. 
\item \idx{\it Smoothly regular}: For each neighborhood $U$ of a 
       point $x\in M$ there exists a smooth function $f:M\to \Bbb R$ with 
       $f(x)=1$ and carrier $f\i(\Bbb R\setminus\{0\})$ contained in 
       $U$; equivalently the initial topology with respect to 
       $C^\infty(M,\Bbb R)$ equals the natural topology.
\item \idx{\it Smoothly real compact}: Any bounded algebra 
       homomorphism $C^\infty(M,\Bbb R)\to \Bbb R$ is given by a 
       point evaluation; equivalently, the natural mapping 
       $M\to \operatorname{Hom}(C^\infty(M,\Bbb R),\Bbb R)$ is 
       surjective.
\endroster

\subhead \nmb.{2.5}. Submanifolds \endsubhead A subset $N$ of a manifold
$M$ is called a \idx{\it submanifold}, if for each $x\in N$ there is
a chart $(U,u)$ of $M$ such that 
$u(U\cap N) = u(U)\cap F_U$, where $F_U$ is a $c^\infty$-closed linear
subspace of the convenient model space $E_U$.
Then clearly $N$ is itself a manifold with $(U\cap
N,u\restriction U\cap N)$ as charts, where $(U,u)$ runs through
all \idx{\it submanifold charts} as above.

A submanifold $N$ of $M$ is called a \idx{\it splitting
submanifold} if there is a cover of $N$ by submanifold charts
$(U,u)$ as above such that the $F_U\subset E_U$ are complemented
(i.e. splitting) linear subspaces. Then obviously every
submanifold chart is splitting.

\subhead \nmb.{2.6}. Products \endsubhead
Let $M$ and $N$ be smooth manifolds
described by smooth atlases $(U_\al,u_\al)_{\al\in A}$ and
$(V_\be,v_\be)_{\be\in B}$, respectively. Then the family
$(U_\al\x V_\be,u_\al\x v_\be:U_\al\x V_\be \to E_\al\x 
F_\be)_{(\al,\be)\in A\x B}$ is a smooth atlas for the 
cartesian product $M\x N$. 
Beware, however, the manifold topology of $M\x N$ may be finer than 
the product topology, see \cit!{10}. If $M$ and $N$ are metrizable, 
then it is the product topology, by \cit!{10} again. 
Clearly the projections 
$$ M @<{pr_1}<< M\x N @>{pr_2}>> N$$
are also smooth. The \idx{\it product\ign{ of manifolds}} $(M\x N, pr_1, pr_2)$ has the
following universal property:

For any smooth manifold $P$ and smooth mappings $f:P\to M$ and
$g:P\to N$ the mapping $(f,g):P\to M\x N$,
$(f,g)(x)=(f(x),g(x))$, is the unique smooth mapping
with $pr_1\o (f,g) = f$, $pr_2\o (f,g) = g$.

\proclaim{\nmb.{2.7}. Lemma} \cit!{10}
For a convenient vector space $E$ and any smooth manifold $M$ the
set $C^\infty(M,E)$ of smooth $E$-valued functions on $M$
is again a convenient vector space with the 
locally convex topology of uniform convergence on compact subsets of 
compositions with smooth curves in $M$, in all derivatives separately.
Moreover, with this structure, for two manifolds $M$, $N$, the 
exponential law holds:
$$C^\infty(M,C^\infty(N,E))\cong C^\infty(M\x N,E).$$
\endproclaim

\head\totoc\nmb0{3}. Weil functors on infinite dimensional manifolds
\endhead

\subheading{\nmb.{3.1}} A real commutative
algebra $A$ with unit 1 is called 
\idx{\it formally real\ign{ commutative algebra}} if for
any $\row a1n\in A$ the element $1+ a_1^2+\dots+a_n^2$ is
invertible in $A$. 
Let $E=\{e\in A: e^2=e,e\ne0\}\subset A$ be the set of
all nonzero idempotent elements in $A$. It is not empty since 
$1\in E$. An idempotent $e\in E$ is said to be {\it minimal} if for 
any $e'\in E$ we have $ee'=e$ or $ee'=0$.

\proclaim{Lemma} Let $A$ be a real commutative algebra with unit
which is formally real and finite dimensional as a real vector space. 

Then there is a decomposition $1=e_1+\dots+e_k$ into all minimal 
idempotents.
Furthermore $A=A_1\oplus \dots \oplus A_k$, where 
$A_i=e_iA=\Bbb R\cdot e_i \oplus N_i$, and $N_i$ is a
nilpotent ideal.
\endproclaim

This is standard, see \cit!{6},~35.1, for a proof. 

\subhead \nmb.{3.2}. Definition \endsubhead
A \idx{\it Weil algebra} $A$ is
a real commutative algebra with unit 
which is of the form $A=\Bbb R\cdot 1\oplus N$, where $N$ is a
finite dimensional ideal of nilpotent elements. 

So by lemma \nmb!{3.1} a formally real and finite dimensional
unital commutative algebra is the direct sum of finitely many
Weil algebras.

\subhead \nmb.{3.3}. \idx{Chart description of Weil functors} 
\endsubhead
Let  $A=\Bbb R\cdot 1\oplus N$ be a Weil algebra. We want to
associate to it a functor $T_A:\Mf \to \Mf$ from the category
$\Mf$ of all smooth manifolds modelled on convenient vector spaces 
into itself. 

\remark{Step 1} If $f\in C^\infty(\Bbb R,\Bbb R)$ and $\la1+n\in \Bbb
R\cdot1\oplus N=A$, we consider the Taylor expansion 
$j^\infty f(\la)(t) = \sum_{j=0}^\infty
\frac{f^{(j)}(\la)}{j!}t^j$ of $f$ at $\la$ and we put 
$$T_A(f)(\la1+n):= f(\la)1+\sum_{j=1}^\infty \frac{f^{(j)}(\la)}{j!}n^j,$$
which is finite sum, since $n$ is nilpotent. Then 
$T_A(f):A\to A$ is smooth and we get 
$T_A(f\o g)=T_A(f)\o T_A(g)$ and $T_A(\Id_{\Bbb R})= \Id_A$.
\endremark

\remark{Step 2} If $f\in C^\infty(\Bbb R,F)$ for a convenient vector 
space $F$ and 
$\la1+n\in \Bbb R\cdot1\oplus N=A$, we consider the Taylor expansion 
$j^\infty f(\la)(t) = \sum_{j=0}^\infty
\frac{f^{(j)}(\la)}{j!}t^j$ of $f$ at $\la$ and we put 
$$
T_A(f)(\la1+n):= 1\otimes f(\la)
     +\sum_{j=1}^\infty n^j\otimes \frac{f^{(j)}(\la)}{j!},
$$
which is finite sum, since $n$ is nilpotent. Then 
$T_A(f):A\to A\otimes F =: T_AF$ is smooth.
\endremark

\remark{Step 3} For $f\in C^\infty(E,F)$, where $E$, $F$ are 
convenient vector spaces, we want to
define the value of $T_A(f)$ at an element of the convenient vector 
space $T_AE = A\otimes E$. Such an element may be uniquely written as
$1\otimes x_1 + \sum_j n_j\otimes x_j$, where $1$ and the $n_j\in N$ 
form a fixed 
finite linear basis of $A$, and where the $x_i\in E$.
Let again 
$j^\infty f(x_1)(y)=\sum_{k\ge0}\frac1{k!}d^k f(x_1)(y^k)$
be the Taylor expansion of $f$ at $x_1\in E$ for 
$y\in E$. Then we put 
$$\multline
T_A(f)(1\otimes x_1 + \sum_j n_j\otimes x_j):=\\
= 1\otimes f(x_1) + \sum_{k\ge0}\frac1{k!}\sum_{j_1,\dots,j_k} 
 n_{j_1}\dots n_{j_k}\otimes d^kf(x_1)(x_{j_1},\dots,x_{j_k})
\endmultline$$
which is again a finite sum. A change of basis in $N$ induces 
the transposed change in the $x_i$, namely 
$\sum_i(\sum_j a_i^j n_j)\otimes \bar x_i = \sum_j n_j 
\otimes (\sum_i a_i^j \bar x_i)$, so the value of $T_A(f)$ is 
independent of the choice of the basis of $N$. 
Since the Taylor expansion of a composition is the composition
of the Taylor expansions we have $T_A(f\o g)= T_A(f)\o T_A(g)$
and $T_A(\Id_{E})= \Id_{T_AE}$. 

If $\ph: A\to B$ is a homomorphism between two Weil algebras we
have $(\ph\otimes F)\o T_Af = T_Bf\o(\ph\otimes E)$ for 
$f\in C^\infty(E,F)$.
\endremark

\remark{Step 4} Let $\pi =\pi _A:A\to A/N= \Bbb R$ 
be the projection onto the
quotient field of the Weil algebra $A$. This is a surjective
algebra homomorphism, so by step 3 the following diagram
commutes for $f\in C^\infty(E,F)$:
\endremark
$$\newCD
A\otimes E  @()\L{T_Af}@(1,0) @()\L{\pi\otimes E}@(0,-1) & 
      A\otimes F @()\l{\pi\otimes F}@(0,-1)\\
E @()\L f@(1,0) & F
\endnewCD$$
If $U\subset E$ is a $c^\infty$-open subset we put
$T_A(U):=(\pi\otimes E)\i(U) = (1\otimes U)\x ( N\otimes E)$, 
which is an $c^\infty$-open subset in 
$T_A(E):= A\otimes E$. If $f:U\to V$ is a smooth mapping between
$c^\infty$-open subsets $U$ and $V$ of $E$ and $F$,
respectively, then the construction of step 3, applied to
the Taylor expansion of $f$ at points in $U$, produces a smooth
mapping $T_Af:T_AU\to T_AV$, which fits into the following
commutative diagram:
$$\cgaps{.6;1.2;.6}\newCD
U\x ( N\otimes E) @()\a=@(1,0) @()\l{pr_1}@(1,-1) & 
	T_AU  @()\L{T_Af}@(1,0) @()\l{\pi\otimes E}@(0,-1) & 
	T_AV @()\a=@(1,0) @()\L{\pi\otimes F}@(0,-1) & V\x ( N\otimes F) 
     @()\l{pr_1}@(-1,-1)\\
& U @()\L f@(1,0) &  V  & 
\endnewCD$$
We have $T_A(f\o g) = T_Af\o T_Ag$ and $T_A(\Id_U) = \Id_{T_AU}$,
so $T_A$ is now a covariant functor on the category of $c^\infty$-open
subsets of convenient vector spaces and smooth mappings between them. 

\remark{Step 5} 
Let $M$ be a smooth manifold, 
let $(U_\al,u_\al:U_\al\to u_\al(U_\al)\subset E_\al)$
be a smooth atlas of $M$ with chart changings 
$u_{\al\be} := u_\al\o u_\be\i:u_\be(U_{\al\be})\to u_\al(U_{\al\be})$. 
Then the smooth mappings
$$\cgaps{1.5}\newCD
T_A(u_\be(U_{\al\be})) @()\L{T_A(u_{\al\be})}@(1,0) 
     @()\L{\pi\otimes E_\be}@(0,-1) & 
	T_A(u_\al(U_{\al\be}))  @()\l{\pi\otimes E_\al}@(0,-1) \\
u_\be(U_{\al\be}) @()\L{u_{\al\be}}@(1,0) & u_\al(U_{\al\be}) 
\endnewCD$$ form again a cocycle of chart changings and we may use 
them to glue the $c^\infty$-open sets $T_A(u_\al(U_\al)) = 
u_\al(U_\al)\x ( N\otimes E_\al) \subset T_AE_\al$
in order to obtain a smooth manifold which we denote by $T_AM$.
By the diagram above we see that $T_AM$ will be the total space
of a fiber bundle $T(\pi_A,M)=\pi_{A,M}:T_AM \to M$, since the atlas
$(T_A(U_\al),T_A(u_\al))$ constructed just now is already a
fiber bundle atlas. So if $M$ is Hausdorff then also $T_AM$ is 
Hausdorff, since two points
$x_i$ can be separated in one chart if they are in the same
fiber, or they can be separated by inverse images under
$\pi _{A,M}$ of open sets in $M$ separating their projections.

This construction does not depend on the choice of the atlas.
For two atlases have a common refinement and one may pass to this.

If $f\in C^\infty(M,M')$ for two manifolds $M$, $M'$, we apply
the functor $T_A$ to the local representatives of $f$ with
respect to suitable atlases. This gives local representatives
which fit together to form a smooth mapping $T_Af:T_AM\to T_AM'$.
Clearly we again have $T_A(f\o g) = T_Af\o T_Ag$ and $T_A(\Id_M)=
\Id_{T_AM}$, so that $T_A:\Mf\to \Mf$ is a covariant functor.
\endremark

\subhead \nmb.{3.4}. Remark \endsubhead If we apply the construction of
\nmb!{3.3}, step 5 to the algebra $A=0$, which we did not
allow ($1\neq 0\in A$), then $T_0M$ depends on the choice of the
atlas. If each chart is connected, then $T_0M = \pi_0(M)$,
computing the connected components of $M$. If each chart meets
each connected component of $M$, then $T_0M$ is one point.

\proclaim{\nmb.{3.5}. Theorem. \idx{Main properties of Weil functors}} 
Let $A = \Bbb R\cdot 1\oplus N$
be a Weil algebra, where $N$ is the maximal ideal of nilpotents.
Then we have:

1. The construction of \nmb!{3.3} defines a covariant functor 
$T_A:\Mf\to \Mf$ such that $(T_AM,\pi_{A,M}, M)$ is a
smooth fiber bundle with standard fiber $ N\otimes E$ if $M$ is 
modelled on the convenient space $E$. For any
$f\in C^\infty(M,M')$ we have a commutative diagram
$$\newCD 
T_AM @()\L{T_Af}@(1,0) @()\L{\pi_{A,M}}@(0,-1) & 
	T_AM'  @()\l{\pi_{A,M'}}@(0,-1) \\
M    @()\L f@(1,0) &  M'\rlap.
\endnewCD$$
So $(T_A, \pi_A)$ is a bundle functor on $\Mf$, which gives a 
vector bundle on $\Mf$ if and only if $N$ is
nilpotent of order 2.

2. The functor $T_A:\Mf\to \Mf$ is multiplicative: it respects
products. It maps the following classes of mappings into itself:
immersions, splitting immersions, embeddings, splitting embeddings, 
closed embeddings, submersions, splitting submersions, 
surjective submersions, fiber bundle projections.
It also respects transversal pullbacks.
For fixed manifolds $M$ and $M'$ the mapping 
$T_A:C^\infty(M,M') \to C^\infty(T_AM,T_AM')$ is smooth, so it
maps smoothly parametrized families to smoothly parametrized
families. 

3. If $(U_\al)$ is an open cover of $M$ then $T_A(U_\al)$ is
also an open cover of $T_AM$.

4. Any algebra homomorphism $\ph:A\to B$ between Weil algebras
induces a natural transformation $T(\ph,\quad)=T_\ph:T_A\to T_B$. 
If $\ph$ is injective, then $T(\ph,M):T_AM\to T_BM$ is a closed
embedding for each manifold $M$. If $\ph$ is surjective, then
$T(\ph,M)$ is a fiber bundle projection for each $M$.
So we may view $T$ as a co-covariant bifunctor from the category
of Weil algebras times $\Mf$ to $\Mf$.
\endproclaim

\demo{Proof} 1. The main assertion is clear from \nmb!{3.3}.
The fiber bundle $\pi_{A,M}:T_AM\to M$ is a vector bundle if and only
if the transition functions $T_A(u_{\al\be})$ are fiber linear
$ N\otimes E_\al\to  N\otimes E_\be$. So only the first derivatives of
$u_{\al\be}$ should act on $N$, so any product of two elements
in $N$ must be 0, thus $N$ has to be nilpotent of order 2.

2. The functor $T_A$ respects finite products in the category of 
$c^\infty$-open
subsets of convenient vector spaces by \nmb!{3.3}, step 3 and 5. All 
the
other assertions follow by looking again at the chart structure
of $T_AM$ and by taking into account that $f$ is part of $T_Af$
(as the base mapping).

3. This is obvious from the chart structure.

4. We define $T(\ph,E) := \ph\otimes E: A\otimes E \to B\otimes E$. 
By
\nmb!{3.3}, step 3, this restricts to a natural transformation
$T_A\to T_B$ on the category of $c^\infty$-open subsets of convenient 
vector spaces and
by gluing also on the category $\Mf$. Obviously $T$ is a
co-covariant bifunctor on the indicated categories. Since
$\pi_B\o\ph=\pi_A$ ($\ph$ respects the identity), we have 
$T(\pi_B,M)\o T(\ph,M) =T(\pi_A,M)$, so $T(\ph,M): T_AM\to T_BM$
is fiber respecting for each manifold $M$. In each fiber chart it is 
a linear mapping
on the typical fiber $ N_A\otimes E\to  N_B\otimes E$. 

So if $\ph$ is injective, $T(\ph,M)$ is fiberwise injective and
linear in each canonical fiber chart, so it is a closed embedding.

If $\ph$ is surjective, let $N_1:= \ker\ph\subseteq N_A$, and
let $V\subset N_A$ be a linear complement to $N_1$. Then if $M$ is 
modeled on convenient vector spaces $E_\al$ 
and for the canonical charts we have the commutative diagram:
$$\rgaps{1;1;.4}\cgaps{2}\newCD
T_AM @()\L{T(\ph,M)}@(1,0) & T_BM \\
T_A(U_\al) @()\L{T(\ph,U_\al)}@(1,0) @(0,1) @()\L{T_A(u_\al)}@(0,-1) &  
	T_B(U_\al) @(0,1)   @()\L{T_B(u_\al)}@(0,-1)\\
u_\al(U_\al)\x ( N_A\otimes E_\al)  @()\L{\Id\x (\ph|N_A\otimes E_\al)}@(1,0) 
     @()\a=@(0,-1) &  
	u_\al(U_\al)\x ( N_B\otimes E_\al) @()\a=@(0,-1)\\
u_\al(U_\al)\x ( N_1\otimes E_\al)\x (V\otimes E_\al)  
     @()\L{\Id\x 0\x Iso}@(1,0) &  
	u_\al(U_\al)\x 0 \x ( N_B\otimes E_\al)
\endnewCD$$
So $T(\ph,M)$ is a fiber bundle projection with standard
fiber $E_\al\otimes \ker \ph$.
\qed\enddemo

\proclaim{\nmb.{3.6}. Theorem} Let $A$ and $B$ be Weil algebras.
Then we have:
\roster
\item We get the algebra $A$ back from the Weil functor $T_A$ by
    $T_A(\Bbb R)=A$ with addition $+_A=T_A(+_{\Bbb R})$,
    multiplication $m_A = T_A(m_{\Bbb R})$ and scalar multiplication
    $m_t = T_A(m_t):A\to A$.
\item The natural transformations $T_A\to T_B$ correspond
    exactly to the algebra homomorphisms $A\to B$%:
\endroster
\endproclaim

\demo{Proof} \therosteritem1 is obvious.
\therosteritem2 For a natural transformation $\ph:T_A\to T_B$
its value $\ph_{\Bbb R}:T_A(\Bbb R)=A\to T_B(\Bbb R)=B$ is an
algebra homomorphisms. The inverse of this mapping is already
described in theorem \nmb!{3.5}.4.
\qed\enddemo

\proclaim{\nmb.{3.7}. Proposition} For two manifolds $M_1$ and
$M_2$, with $M_2$ smoothly real compact and smoothly regular, the mapping 
$$\gather C^\infty(M_1,M_2)\to
    \operatorname{Hom}(C^\infty(M_2,\Bbb R),C^\infty(M_1,\Bbb R))\\
f\mapsto (f^*:g\mapsto g\o f)
\endgather$$
is bijective.
\endproclaim

\demo{Proof} Let $x_1\in M_1$ and
$\ph\in\operatorname{Hom}(C^\infty(M_2,\Bbb R),C^\infty(M_1,\Bbb R))$.
Then $ev_{x_1}\o \ph$ is in
$\operatorname{Hom}(C^\infty(M_2,\Bbb R),\Bbb R)$, so by
\nmb!{2.4} there is a $x_2\in M_2$ such that 
$ev_{x_1}\o \ph = ev_{x_2}$ since $M_2$ is smoothly real compact, and 
$x_2$ is unique since $M_2$ is smoothly Hausdorff. 
If we write $x_2 = f(x_1)$, then
$f:M_1\to M_2$ and $\ph(g)= g\o f$ for all $g\in
C^\infty(M_2,\Bbb R)$. This implies that $f$ is smooth, since $M_2$ 
is smoothly regular.
\qed\enddemo

\subhead{\nmb.{3.8}. Remark}\endsubhead 
If $M$ is a smoothly real compact and smoothly regular manifold 
we consider the set  
$D_A(M):=\operatorname{Hom}(C^\infty(M,\Bbb R),A)$ of all bounded 
homomorphisms from the convenient algebra of smooth functions on $M$ into 
a Weil algebra $A$. 
Obviously we have a natural mapping $T_AM\to D_AM$ which is given by 
$X\mapsto (f\mapsto T_A(f).X)$, using \nmb!{3.5} and \nmb!{3.6}.

Let $\Bbb D$ be the algebra of Study numbers $\Bbb R.1 \oplus \Bbb R.\de$ 
with $\de^2=0$.  
Then $T_{\Bbb D}M= TM$, the tangent bundle, and 
$D_{\Bbb D}(M)$ is the smooth bundle of all 
operational tangent vectors, i.e\. bounded derivations at a point $x$ of 
the algebra of germs $C^\infty_x(M,\Bbb R)$ see \cit!{10}. 
We want to point out that even on Hilbert spaces there exist 
derivations which are differential operators of order 2 and 3, 
respectively, see \cit!{10}.

It would be nice if $D_A(M)$ were a smooth manifold, not only for 
$A=\Bbb D$. We do not know whether this is true.
The obvious method of proof hits severe obstacles, which we now 
explain.  

Let $A=\Bbb R.1\oplus N$ for a nilpotent finite dimensional ideal 
$N$, let $\pi:A\to \Bbb R$ be the corresponding projection.
Then for $\ph\in D_A(M)=\operatorname{Hom}(C^\infty(M,\Bbb R),A)$ the 
character $\pi\o \ph =\ev_x$ for a unique $x\in M$, since $M$ is 
smoothly real compact. Moreover $X:= 
\ph-\ev_x.1:C^\infty(M,\Bbb R)\to N$ satisfies the \idx{\it expansion 
property} at $x$:
$$
X(fg) = X(f).g(x) + f(x).X(g) + X(f).X(g).\tag1
$$
Conversely a bounded linear mapping $X:C^\infty(M,\Bbb R)\to N$ with 
property \thetag1 is called an \idx{\it expansion at $x$}. 
Clearly each expansion at $x$ defines a bounded homomorphism $\ph$ with 
$\pi\o \ph=\ev_x$. So we view $D_A(M)_x$ as the set of all expansions 
at $x$. Note first that for an expansion $X\in D_A(M)_x$ the value 
$X(f)$ depends only on the germ of $f$ at $x$: If $f|U=0$ for a 
neighborhood $U$ of $x$, choose a smooth function $h$ with $h=1$ off 
$U$ and $h(x)=0$. Then $h^kf=f$ and 
$X(f)=X(h^kf)=0+0+X(h^k)X(f)=\dots=X(h)^kX(f)$ which is 0 for $k$ 
larger than the nilpotence index of $N$.

Suppose now that $M=U$ is a $c^\infty$-open subset of a convenient 
vector space $E$. 
We can ask whether $D_A(U)_x$ is a smooth manifold. 
We have no proof of this. Let us sketch the difficulty.
A natural way to prove that would be 
by induction on the nilpotence index of $N$. Let
$N_0:=\{n\in N: n.N=0\}$, which is an ideal in $A$. 
Consider the short exact sequence
$$
0\to N_0\to N @>{p}>> N/N_0\to 0
$$
and a linear section $s:N/N_0\to N$. 
For $X:C^\infty(U,\Bbb R)\to N$ we consider $\bar X:= p\o X$ and 
$X_0:=X-s\o\bar X$. Then $X$ is an expansion at $x\in U$ if and only if 
\roster
\item""
$\bar X$ is an expansion at $x$ with values in $N/N_0$ and $X_0$ 
satisfies 
$$X_0(fg) = X_0(f)g(x)+f(x)X_0(g) + s(\bar X(f)).s(\bar X(g)) 
     - s(\bar X(f).\bar X(g)).\tag 2$$
\endroster
Note that \thetag2 is an affine equation in $X_0$ for fixed $\bar X$. 
By induction the $\bar X\in D_{A/N_0}(U)_x$ 
form a smooth manifold, and the fiber over a 
fixed $\bar X$ consists of all $X=X_0+s\o \bar X$ with $X_0$ 
in the closed affine subspace described by \thetag2, whose model 
vector space is the space of all derivations at $x$. 
If we were able to find a (local) section $D_{A/N_0}(U)\to D_A(U)$ 
and if these sections would fit together nicely 
we could then conclude that $D_A(U)$ were the total space of a smooth 
affine bundle over $D_{A/N_0}(U)$, so it would be smooth. 
But this translates to a lifting problem as follows:
A homomorphism $C^\infty(U,\Bbb R)\to A/N_0$ has to be lifted in a 
`natural way' to $C^\infty(U,\Bbb R)\to A$. But we know that in 
general $C^\infty(U,\Bbb R)$ is not a free $C^\infty$-algebra, see 
\nmb!{4.4} for comparison. 

%A diffeomorphism $\ps:U\to U'$ induces an isomorphism of iterated affine 
%bundles $D_A(\ps)_x:D_A(U)_x \to D_A(U')_{\ps(x)}$ by 
%$(D_A(\ps)_x.X)(f)=X(\ps^*f)$. Thus $D_A(U)\cong U\x D_A(U)_{x}$ is a 
%smooth manifold, and if $M$ is a smooth manifold with atlas $U_\al$, 
%we may glue the 
%iterated affine bundle $D_A(M)\to M$ from the local parts 
%$D_A(U_\al)$, as for the operational tangent bundles.

\subhead{\nmb.{3.9}}\endsubhead
The basic facts from the theory of Weil 
functors are completed by the following assertion.

\proclaim{Proposition} Given two Weil algebras $A$ and $B$, the 
composed functor $T_A\o T_B$ is a Weil functor generated by the 
tensor product $A\otimes B$.
\endproclaim

\demo{Proof}
For a convenient vector space $E$ we have $T_A(T_BE)= 
A\otimes B\otimes E$ and this is compatible with the action of smooth 
mappings, by \nmb!{3.3}.
\qed\enddemo

\proclaim{Corollary} There is a canonical 
natural equivalence $T_A\o T_B\cong T_B\o T_A$ generated by the 
exchange algebra isomorphism $A\otimes B\cong B\otimes A$.
\endproclaim

\subhead{\nmb.{3.10}. Weil functors and Lie groups}\endsubhead
We shall use the notion of a regular infinite dimensional Lie group, 
modelled on convenient vector spaces, as laid out in \cit!{9}, 
following the lead of Omori et.\ al.\ \cit!{20} and Milnor \cit!{15}. 
We just remark that they have unique smooth exponential mappings, and 
that no smooth Lie group is known which is not regular. 
We shall use the notation $\mu:G\x G\to G$ for the multiplication and 
$\nu:G\to G$ for the inversion. 
The tangent bundle $TG$ of a regular Lie 
group $G$ is again a Lie group, the semidirect product $\frak g 
\ltimes G$ of $G$ with its Lie algebra $\frak g$.

Now let $A$ be a Weil algebra and let $T_A$ be its Weil functor. Then 
the space 
$T_A(G)$ is again a Lie group with multiplication $T_A(\mu)$ and inversion
$T_A(\nu)$. By the properties \nmb!{3.5} of the Weil functor $T_A$ 
we have a surjective homomorphism $\pi_A:T_AG\to G$ of Lie groups. 
Following the analogy with the tangent bundle, for $a\in G$ we will 
denote its fiber over $a$ by $(T_A)_aG\subset T_AG$, likewise for 
mappings. With this notation we have the following commutative 
diagram, where we assume that $G$ is a regular Lie group: 
$$\cgaps{0.5;1;1;0.5}\rgaps{0.7;1}\newCD
& \frak g \otimes N @(1,0) @()\a=@(0,-1) & 
	\frak g\otimes A  @()\a=@(0,-1) & & \\
0 @(1,0) & (T_A)_0\frak g @(1,0) @()\L{(T_A)_0\exp^G}@(0,-1) & 
	T_A\frak g @(1,0) @()\L{T_A\exp^G}@(0,-1) & 
	\frak g @(1,0) @()\l{\exp^G}@(0,-1) & 0\\
e @(1,0) & (T_A)_eG @(1,0) & T_AG @()\L{\pi_A}@(1,0) & G @(1,0) & e
\endnewCD$$
The structural mappings (Lie bracket, exponential mapping, 
evolution operator, adjoint action) are determined by multiplication 
and inversion. Thus their images under the Weil functor $T_A$ are again 
the same structural mappings. But note 
that the canonical flip mappings have to be inserted like follows.    
So for example 
$$\frak g\otimes A \cong T_A\frak g = T_A(T_eG) @>{\ka}>> T_e(T_AG)$$ 
is the Lie algebra of 
$T_AG$ and the Lie bracket is just $T_A([\quad,\quad])$. 
Since the bracket is bilinear, the description of \nmb!{3.3} implies 
that $[X\otimes a,Y\otimes b]_{T_A\frak g}=[X,Y]_{\frak g}\otimes ab$. 
Also $T_A\exp^G = \exp^{T_AG}$. If $\exp^G$ is a diffeomorphism near 
0, $(T_A)_0(\exp^G):(T_A)_0\frak g\to (T_A)_eG$ is also a diffeomorphism 
near 0, since $T_A$ is local.
The natural transformation $0_G:G\to T_AG$ is a homomorphism 
which splits the bottom row of the diagram, so $T_AG$ is the 
semidirect product $(T_A)_0\frak g\ltimes G$ via the mapping
$T_A\rho:(u,g)\mapsto T_A(\rho_g)(u)$. So from \cit!{9}, theorem~5.5, 
we may 
conclude that $T_AG$ is again a regular Lie group, if $G$ is regular.   
If $\om^G:TG\to T_eG$ is the Maurer Cartan form of $G$ (i\.e\. the left 
logarithmic derivative of $\Id_G$) then 
$$\ka_0\o T_A\om^G\o \ka:TT_AG\cong T_ATG\to T_AT_eG\cong T_eT_AG$$
is the Maurer Cartan form of $T_AG$.

\head\totoc\nmb0{4}. Product preserving functors from finite 
dimensional manifolds to infinite dimensional ones \endhead

\subhead\nmb.{4.1}. Product preserving functors \endsubhead
Let $\Mf_{\text{fin}}$ denote the category of all finite dimensional 
separable Hausdorff smooth manifolds, with smooth mappings as 
morphisms. 
Let $F:\Mf_{\text{fin}}\to \Mf$ be a functor which preserves products 
in the following sense:
The diagram
$$F(M_1) @<{F(pr_1)}<< F(M_1\x M_2) @>{F(pr_2)}>> F(M_2)$$
is always a product diagram. 

Then $F(\text{point})=\text{point}$, by 
the following argument:
$$\newCD
F(\text{point}) & 
F(\text{point}\x\text{point}) @()\L{F(pr_1)}\l{\cong}@(-1,0) 
	@()\L{F(pr_2)}\l{\cong}@(1,0) & 
F(\text{point})\\
& \text{point} @()\l{f_1}@(-1,1) @()\L{f}@(0,1) @()\l{f_2}@(1,1) 
\endnewCD$$
Each of $f_1$, $f$, and $f_2$ determines each other uniquely, thus 
there is only one mapping $f_1:\text{point}\to F(\text{point})$, so 
the space $F(\text{point})$ is single pointed.

We also require that $F$ has the following two properties:
\roster
\item  The map on morphisms
       $F:C^\infty(\Bbb R^n,\Bbb R)\to C^\infty(F(\Bbb R^n),F(\Bbb R))$ 
       is smooth, where we regard $C^\infty(F(\Bbb R^n),F(\Bbb R))$ 
       as smooth space, see \cit!{2} or \cit!{10}. Equivalently the 
       associated mapping 
       $C^\infty(\Bbb R^n,\Bbb R)\x F(\Bbb R^n)\to F(\Bbb R)$ is  
       smooth.
\item  There is a natural transformation $\pi:F\to \Id$ such that for 
       each $M$ the mapping $\pi_M:F(M)\to M$ is a fiber bundle, and 
       for an open submanifold $U\subset M$ the mapping 
       $F(\operatorname{incl}):F(U)\to F(M)$ is a pullback. 
\endroster

\subhead\nmb.{4.2}. $C^\infty$-algebras \endsubhead
An $\Bbb R$-algebra is a commutative ring $A$ with unit together with 
a ring homomorphism $\Bbb R\to A$. Then every map 
$p:\Bbb R^n\to \Bbb R^m$ which is given by an $m$-tuple of real 
polynomials $(p_1,\dots,p_m)$ can be interpreted as a mapping 
$A(p):A^n\to A^m$ in such a way that projections, composition, and 
identity are preserved, by just evaluating each polynomial $p_i$ on 
an $n$-tuple $(a_1,\dots,a_n)\in A^n$.

A $C^\infty$-algebra $A$ is a real algebra in which we can moreover 
interpret all smooth mappings $f:\Bbb R^n\to \Bbb R^m$. There is a 
corresponding map $A(f):A^n\to A^m$, and again projections, 
composition, and the identity mapping are preserved. 

More precisely, a $C^\infty$-algebra $A$ is a product preserving 
functor from the category $C^\infty$ to the category of sets, where 
$C^\infty$ has as objects all spaces $\Bbb R^n$, $n\ge 0$, and all 
smooth mappings between them as arrows. Morphisms between 
$C^\infty$-algebras are then natural transformations: they correspond 
to those algebra homomorphisms which preserve the interpretation of 
smooth mappings.

Let us explain how one gets the algebra structure from this 
interpretation. Since $A$ is product preserving, we have
$A(\text{point})=\text{point}$. All the laws for a commutative
ring with unit can be formulated by commutative diagrams of
mappings between products of the ring and the point. We do this
for the ring $\Bbb R$ and apply the product preserving functor
$A$ to all these diagrams, so we get the laws for the
commutative ring $A(\Bbb R)$ with unit $A(1)$ with the exception of 
$A(0)\ne A(1)$ which we will check later for the case 
$A(\Bbb R)\ne\text{point}$. 
Addition is given by $A(+)$
and multiplication by $A(m)$.
For $\la\in \Bbb R$ the mapping $A(m_\la):A(\Bbb R)\to A(\Bbb R)$ 
equals multiplication with the element $A(\la)\in A(\Bbb R)$,
since the following diagram commutes:
$$\cgaps{1.7;1}\rgaps{.6;.6;.6}\newCD
A(\Bbb R) @()\L{\cong}@(0,-1) @()\L{A(m_\la)}@(2,-1) & & \\
A(\Bbb R)\x\text{point}  @()\L{\Id\x A(\la)}@(1,0) @()\L{\cong}@(0,-1) & 
	A(\Bbb R)\x A(\Bbb R) @(1,0) @()\a=@(0,-1) & 
     A(\Bbb R) \\
A(\Bbb R\x\text{point}) @()\L{A(\Id\x\la)}@(1,0) & 
	A(\Bbb R\x\Bbb R) @()\l{A(m)}@(1,1) & 
\endnewCD$$
We may investigate now the difference between 
$A(\Bbb R)=\text{point}$ and $A(\Bbb R)\ne \text{point}$. In the 
latter case for $\la\ne0$ we have $A(\la)\ne A(0)$ since 
multiplication by $A(\la)$ equals $A(m_{\la})$ which is a 
diffeomorphism for $\la\ne0$ and factors over a one pointed space for 
$\la=0$. So for $A(\Bbb R)\ne\text{point}$ which we assume from now 
on, the group homomorphism $\la\mapsto A(\la)$ from $\Bbb R$ into 
$A(\Bbb R)$ is actually injective.

This definition of $C^\infty$-algebras is due to Lawvere \cit!{12}, 
for a thorough account see Moerdijk-Reyes \cit!{16}, for a discussion 
from the point of view of functional analysis see \cit!{3}. 
In particular there on a $C^\infty$-algebra $A$ the natural topology 
is defined as the finest locally convex topology on $A$ such that for 
all $a=(a_1,\dots,a_n)\in A^n$ the evaluation mappings 
$\ep_{a}:C^\infty(\Bbb R^n,\Bbb R)\to A$ are continuous.
In \cit!{3}, 6.6 one finds a method to recognize $C^\infty$-algebras 
among locally-m-convex algebras. In \cit!{14} one finds a 
characterization of the algebras of smooth functions on finite 
dimensional algebras among all $C^\infty$-algebras.
 
\proclaim{\nmb.{4.3}. Theorem}
Let $F:\Mf_{\text{fin}}\to \Mf$ be a product preserving 
functor. 
Then either $F(\Bbb R)$ is a point or $F(\Bbb R)$ is a 
$C^\infty$-algebra.
If $\ph:F_1\to F_2$ is a natural transformation between two such 
functors, then 
$\ph_{\Bbb R}:F_1(\Bbb R)\to F_2(\Bbb R)$ is an algebra homomorphism.

If $F$ has property \therosteritem1 then the natural topology on 
$F(\Bbb R)$ is finer than the given manifold 
topology and thus is Hausdorff if the latter is it. 

If $F$ has property \therosteritem2 then $F(\Bbb R)$ is a 
local algebra with an algebra homomorphism 
$\pi=\pi_{\Bbb R}:F(\Bbb R)\to \Bbb R$ whose kernel is the maximal 
ideal. 
\endproclaim

\demo{Proof}
By definition $F$ restricts to a product preserving functor from the 
category of all $\Bbb R^n{}$'s and smooth mappings between them, thus 
it is a $C^\infty$-algebra. 

If $F$ has property \therosteritem1 then for all 
$a=(a_1,\dots,a_n)\in F(\Bbb R)^n$ the evaluation mappings are given 
by
$$
\ep_{a}=\ev_a\o F: 
C^\infty(\Bbb R^n,\Bbb R)\to C^\infty(F(\Bbb R)^n,F(\Bbb R))\to F(\Bbb R)
$$ and thus are even smooth.

If $F$ has property \therosteritem2 then obviously 
$\pi_{\Bbb R}=\pi:F(\Bbb R)\to \Bbb R$ is an algebra homomorphism. It 
remains to show that the kernel of $\pi$ is the largest ideal. So if 
$a\in A$ has $\pi(a)\ne0\in \Bbb R$ then we have to show that $a$ is 
invertible in $A$. 
Since the following diagram is a pullback,
$$\CD
F(\Bbb R\setminus \{0\}) @>{F(i)}>> F(\Bbb R)\\
@V{\pi}VV                            @V{\pi}VV\\
\Bbb R\setminus\{0\} @>i>>           \Bbb R  
\endCD$$
we may assume that $a=F(i)(b)$ for a unique 
$b\in F(\Bbb R\setminus\{0\})$. But then 
$1/i:\Bbb R\setminus\{0\}\to \Bbb R$ is smooth, and $F(1/i)(b)=a\i$, 
since $F(1/i)(b).a=F(1/i)(b).F(i)(b)=F(m)F(1/i,i)(b)=F(1)(b)=1$, 
compare \nmb!{4.2}. 
\qed\enddemo

\subhead\nmb.{4.4}. Examples \endsubhead
Let $A$ be an augmented local $C^\infty$-algebra with maximal ideal 
$N$. Then $A$ is 
quotient of a free $C^\infty$-algebra 
$C^\infty_{\text{fin}}(\Bbb R^\La)$ of smooth functions on some large 
product $\Bbb R^\La$, which depend globally only on finitely many 
coordinates, see \cit!{16} or \cit!{3}. So we have a short exact 
sequence 
$$0\to I \to C^\infty_{\text{fin}}(\Bbb R^\La)@>\ph>> A\to 0.$$
Then $I$ is contained in the codimension 1 maximal ideal 
$\ph\i(N)$, which is easily seen to be
$\{f\in C^\infty_{\text{fin}}(\Bbb R^\la): f(x_0)=0\}$ for some 
$x_0\in \Bbb R^\La$. Then clearly $\ph$ factors over the quotient of 
germs at $x_0$. If $A$ has Hausdorff natural topology, then $\ph$ 
even factors over the Taylor expansion mapping, by the argument in 
\cit!{3},~6.1, as follows. 
Namely, let $f\in C^\infty_{\text{fin}}(\Bbb R^\La)$ 
be infinitely flat at $x_0$. We shall show that $f$ is 
in the closure of the set of all functions with germ 0 at $x_0$.
Let $x_0=0$ without loss. 
Note first that $f$ factors over some quotient 
$\Bbb R^\La\to\Bbb R^N$, and we may replace $\Bbb R^\La$ by 
$\Bbb R^N$ without loss. Define $g:\Bbb R^N\x \Bbb R^N\to \Bbb R^N$,
$$
g(x,y) = \cases 0 &\text{ if }|x|\le|y|,\\
                (1-|y|/|x|)x &\text{ if }|x|>|y|.\endcases
$$
Since $f$ is flat at 0, the mapping $y\mapsto (x\mapsto f_y(x):=f(g(x,y)) $ 
is a continuous mapping $\Bbb R^N\to C^\infty(\Bbb R^N,\Bbb R)$ with 
the property that $f_0=f$ and $f_y$ has germ 0 at 0 for all $y\ne0$.

Thus the augmented local $C^\infty$-algebras whose natural topology 
is Hausdorff are exactly the quotients of algebras of Taylor series 
at 0 of functions in $C^\infty_{\text{fin}}(\Bbb R^\La)$.

It seems that local implies augmented: one has to show that a 
$C^\infty$-algebra which is a field is 1-dimensional. Is this true?

\subhead \nmb.{4.5}. \idx{Chart description of functors induced by 
$C^\infty$-algebras} 
\endsubhead
Let  $A=\Bbb R\cdot 1\oplus N$ be an augmented local $C^\infty$-algebra 
which carries a compatible convenient structure, i.e\. $A$ is a 
convenient vector space and each mapping 
$A:C^\infty(\Bbb R^n,\Bbb R^m)\to C^\infty(A^n,A^m)$ is a well 
defined smooth mapping.
As in the proof of \nmb!{4.3} one sees that the natural topology on 
$A$ is then finer than the given convenient one, thus is Hausdorff. 
Let us call this an \idx{\it augmented local convenient $C^\infty$-algebra}. 

We want to associate to $A$ a functor $T_A:\Mf_{\text{fin}} \to \Mf$ 
from the category $\Mf_{\text{fin}}$ of all finite dimensional 
separable smooth manifolds to the category of smooth manifolds 
modelled on convenient vector spaces.

\remark{Step 1} Let $\pi =\pi _A:A\to A/N= \Bbb R$ 
be the augmentation mapping. This is a surjective
homomorphism of $C^\infty$-algebras, so the following diagram
commutes for $f\in C^\infty(\Bbb R^n,\Bbb R^m)$:
\endremark
$$\newCD
A^n @()\L{T_Af}@(1,0) @()\L{\pi^n}@(0,-1) & 
      A^m @()\l{\pi^m}@(0,-1)\\
\Bbb R^n @()\L f@(1,0) & \Bbb R^m
\endnewCD$$
If $U\subset \Bbb R^n$ is an open subset we put
$T_A(U):=(\pi^n)\i(U) = U\x N^n$, 
which is open subset in 
$T_A(\Bbb R^n):= A^n$. 

\remark{Step 2}
Now suppose that $f:\Bbb R^n\to \Bbb R^m$ 
vanishes on some open set $V\subset \Bbb R^n$.
We claim that then $T_Af$ vanishes on the open set 
$T_A(V)=(\pi^n)\i(V)$. To see this let $x\in V$, and choose a smooth 
function $g\in C^\infty(\Bbb R^n,\Bbb R)$ with $g(x)=1$ and support 
in $V$. Then $g.f=0$ thus 
we have also $0= A(g.f) = A(m)\o A(g,f) = A(g).A(f)$, where the last 
multiplication is pointwise diagonal multiplication between $A$ and 
$A^m$. For $a\in A^n$ with $(\pi^n)(a)=x$ we get 
$\pi(A(g)(a))=g(\pi(a))=g(x)=1$, thus $A(g)(a)$ is invertible in the 
algebra $A$, and from $A(g)(a).A(f)(a)=0$ we may conclude that 
$A(f)(a)=0\in A^m$. 
\endremark

\remark{Step 3}
Now let $f:U\to W$ be a smooth mapping between open sets 
$U\subseteq \Bbb R^n$ and $W\subseteq \Bbb R^m$ . Then we can define 
$T_A(f):T_A(U)\to T_A(W)$ in the following way. 
For $x\in U$ let $f_x:\Bbb R^n\to \Bbb R^m$ be a smooth mapping which 
coincides with $f$ in a neighborhood $V$ of $x$ in $U$. Then by 
step~2 the restriction of $A(f_x)$ to $T_A(V)$ does not depend on the 
choice of the extension $f_x$, and by a standard argument one can 
define uniquely a smooth mapping $T_A(f):T_A(U)\to T_A(V)$. Clearly 
this gives us an extension of the functor $A$ from the category of 
all $\Bbb R^n$'s and smooth mappings into convenient vector spaces to 
a functor from open subsets of $\Bbb R^n$'s and smooth mappings into 
the category of $c^\infty$-open (indeed open) subsets of convenient 
vector spaces. 
\endremark

\remark{Step 4} 
Let $M$ be a smooth finite dimensional manifold, 
let $(U_\al,u_\al:U_\al\to u_\al(U_\al)\subset \Bbb R^m)$
be a smooth atlas of $M$ with chart changings 
$u_{\al\be} := u_\al\o u_\be\i:u_\be(U_{\al\be})\to u_\al(U_{\al\be})$. 
Then by step 3 we get smooth mappings between $c^\infty$-open subsets 
of convenient vector spaces
$$\cgaps{1.5}\newCD
T_A(u_\be(U_{\al\be})) @()\L{T_A(u_{\al\be})}@(1,0) 
     @()\L{\pi}@(0,-1) & 
	T_A(u_\al(U_{\al\be}))  @()\l{\pi}@(0,-1) \\
u_\be(U_{\al\be}) @()\L{u_{\al\be}}@(1,0) & u_\al(U_{\al\be}) 
\endnewCD$$
form again a cocycle of chart changings and we may use 
them to glue the $c^\infty$-open sets $T_A(u_\al(U_\al)) = 
\pi_{\Bbb R^m}\i(u_\al(U_\al)) \subset A^m$
in order to obtain a smooth manifold which we denote by $T_AM$.
By the diagram above we see that $T_AM$ will be the total space
of a fiber bundle $T(\pi_A,M)=\pi_{A,M}:T_AM \to M$, since the atlas
$(T_A(U_\al),T_A(u_\al))$ constructed just now is already a
fiber bundle atlas. So if $M$ is Hausdorff then also $T_AM$ is 
Hausdorff, since two points
$x_i$ can be separated in one chart if they are in the same
fiber, or they can be separated by inverse images under
$\pi _{A,M}$ of open sets in $M$ separating their projections.

This construction does not depend on the choice of the atlas.
For two atlases have a common refinement and one may pass to this.

If $f\in C^\infty(M,M')$ for two manifolds $M$, $M'$, we apply
the functor $T_A$ to the local representatives of $f$ with
respect to suitable atlases. This gives local representatives
which fit together to form a smooth mapping $T_Af:T_AM\to T_AM'$.
Clearly we again have $T_A(f\o g) = T_Af\o T_Ag$ and $T_A(\Id_M)=
\Id_{T_AM}$, so that $T_A:\Mf\to \Mf$ is a covariant functor.
\endremark

\proclaim{\nmb.{4.6}. Theorem. \idx{Main properties}} 
Let  $A=\Bbb R\cdot 1\oplus N$ be a local augmented convenient 
$C^\infty$-algebra.  
Then we have:

1. The construction of \nmb!{4.5} defines a covariant functor 
$T_A:\Mf_{\text{fin}}\to \Mf$ such that $(T_AM,\pi_{A,M}, M)$ is a
smooth fiber bundle with standard fiber $ N^m$ if $\dim M = m$. 
For any $f\in C^\infty(M,M')$ we have a commutative diagram
$$\newCD 
T_AM @()\L{T_Af}@(1,0) @()\L{\pi_{A,M}}@(0,-1) & 
	T_AM'  @()\l{\pi_{A,M'}}@(0,-1) \\
M    @()\L f@(1,0) &  M'\rlap.
\endnewCD$$

2. The functor $T_A:\Mf\to \Mf$ is multiplicative: it respects
products. It respects immersions, embeddings, etc, similarly as in 
\nmb!{3.5}. It also respects transversal pullbacks.
For fixed manifolds $M$ and $M'$ the mapping 
$T_A:C^\infty(M,M') \to C^\infty(T_AM,T_AM')$ is smooth.

3. Any bounded algebra homomorphism $\ph:A\to B$ between augmented 
convenient $C^\infty$-algebras
induces a natural transformation $T(\ph,\quad)=T_\ph:T_A\to T_B$. 
If $\ph$ is split injective, then $T(\ph,M):T_AM\to T_BM$ is a 
split embedding for each manifold $M$. If $\ph$ is split surjective, then
$T(\ph,M)$ is a fiber bundle projection for each $M$.
So we may view $T$ as a co-covariant bifunctor from the category
of augmented convenient $C^\infty$-algebras algebras times 
$\Mf_{\text{fin}}$ to $\Mf$.
\endproclaim

\demo{Proof} 1. The main assertion is clear from \nmb!{4.5}.
The fiber bundle $\pi_{A,M}:T_AM\to M$ is a vector bundle if and only
if the transition functions $T_A(u_{\al\be})$ are fiber linear
$ N\otimes E_\al\to  N\otimes E_\be$. So only the first derivatives of
$u_{\al\be}$ should act on $N$, so any product of two elements
in $N$ must be 0, thus $N$ has to be nilpotent of order 2.

2. The functor $T_A$ respects finite products in the category of 
$c^\infty$-open
subsets of convenient vector spaces by \nmb!{3.3}, step 3 and 5. All 
the
other assertions follow by looking again at the chart structure
of $T_AM$ and by taking into account that $f$ is part of $T_Af$
(as the base mapping).

3. We define $T(\ph,\Bbb R^n) := \ph^n: A^n \to B^n$. 
By \nmb!{4.5}, step 3, this restricts to a natural transformation
$T_A\to T_B$ on the category of open subsets of $\Bbb R^n$'s
by gluing also on the category $\Mf$. Obviously $T$ is a
co-covariant bifunctor on the indicated categories. Since
$\pi_B\o\ph=\pi_A$ ($\ph$ respects the identity), we have 
$T(\pi_B,M)\o T(\ph,M) =T(\pi_A,M)$, so $T(\ph,M): T_AM\to T_BM$
is fiber respecting for each manifold $M$. In each fiber chart it is 
a linear mapping
on the typical fiber $ N_A^m\to  N_B^m$. 

So if $\ph$ is split injective, $T(\ph,M)$ is fiberwise split injective and
linear in each canonical fiber chart, so it is a split embedding.

If $\ph$ is split surjective, let $N_1:= \ker\ph\subseteq N_A$, and
let $V\subset N_A$ be a topological linear complement to $N_1$. Then for
$m=\dim M$ and for the canonical charts we have the commutative diagram:
$$\rgaps{1;1;.4}\cgaps{2}\newCD
T_AM @()\L{T(\ph,M)}@(1,0) & T_BM \\
T_A(U_\al) @()\L{T(\ph,U_\al)}@(1,0) @(0,1) @()\L{T_A(u_\al)}@(0,-1) &  
	T_B(U_\al) @(0,1)   @()\L{T_B(u_\al)}@(0,-1)\\
u_\al(U_\al)\x N_A^m  @()\L{\Id\x \ph|N_A^m}@(1,0) 
     @()\a=@(0,-1) &  
	u_\al(U_\al)\x N_B^m @()\a=@(0,-1)\\
u_\al(U_\al)\x N_1^m\x V^m  @()\L{\Id\x 0\x Iso}@(1,0) &  
	u_\al(U_\al)\x 0 \x N_B^m
\endnewCD$$
So $T(\ph,M)$ is a fiber bundle projection with standard
fiber $E_\al\otimes \ker \ph$.
\qed\enddemo

\proclaim{\nmb.{4.7}. Theorem} Let $A$ and $B$ be augmented 
convenient $C^\infty$-algebras.
Then we have:
\roster
\item We get the convenient $C^\infty$-algebra $A$ back from the 
    functor $T_A$ by restricting to the subcategory of $\Bbb R^n$'s.
\item The natural transformations $T_A\to T_B$ correspond
    exactly to the bounded $C^\infty$-algebra homomorphisms 
    $A\to B$. 
\endroster
\endproclaim

\demo{Proof} \therosteritem1 is obvious.
\therosteritem2 For a natural transformation $\ph:T_A\to T_B$ (which 
is smooth)
its value $\ph_{\Bbb R}:T_A(\Bbb R)=A\to T_B(\Bbb R)=B$ is a 
$C^\infty$-algebra homomorphism which is smooth and thus bounded. 
The inverse of this mapping is already
described in theorem \nmb!{4.6}.4.
\qed\enddemo

\proclaim{\nmb.{4.8}. Proposition}
Let  $A=\Bbb R\cdot 1\oplus N$ be a local augmented convenient 
$C^\infty$-algebra and let $M$ be a smooth finite dimensional 
manifold.  

Then there exists a bijection
$$\ep:T_A(M)\to \Hom(C^\infty(M,\Bbb R),A)$$
onto the space of bounded algebra homomorphisms, 
which is natural in $A$ and $M$. Via $\ep$ the expression 
$\Hom(C^\infty(\quad,\Bbb R),A)$ describes the functor $T_A$ in a 
coordinate free manner.
\endproclaim

\demo{Proof}
{\bf Step 1.} Let $M=\Bbb R^n$, so $T_A(\Bbb R^n)=A^n$. Then for 
$a=(a_1,\dots,a_n)\in A^n$ we have $\ep(a)(f)=A(f)(a_1,\dots,a_n)$, 
which gives a bounded algebra homomorphism $C^\infty(\Bbb R^n,\Bbb R)\to A$.
Conversely, for $\ph\in\Hom(C^\infty(\Bbb R^n,\Bbb R),A)$ consider 
$a=(\ph(\pr_1),\dots,\ph(\pr_n))\in A^n$. Since polynomials are dense 
in $C^\infty(\Bbb R^n,\Bbb R)$, $\ph$ is bounded, and $A$ is 
Hausdorff, $\ph$ is uniquely determined by its values on the
coordinate functions $\pr_i$ (compare \cit!{3},~2.4.(3)), thus 
$\ph(f)=A(f)(a)$ and $\ep$ is bijective. 
Obviously $\ep$ is natural in $A$ and $\Bbb R^n$. 

\smallskip
\noindent{\bf Step 2.} Now let $i:U\subset \Bbb R^n$ be an embedding of 
an open subset. 
Then the image of the mapping 
$$\operatorname{Hom}(C^\infty(U,\Bbb R),A) @>{(i^*)^*}>>
\operatorname{Hom}(C^\infty(\Bbb R^n,\Bbb R),A) 
     @>{\ep_{\Bbb R^n,A}\i}>> A^n$$
is the set $\pi_{A,\Bbb R^n}\i(U) = T_A(U) \subset A^n$, and
$(i^*)^*$ is injective.

To see this let $\ph\in\operatorname{Hom}(C^\infty(U,\Bbb R),A)$.
Then $\ph\i(N)$ is the maximal ideal in $C^\infty(U,\Bbb R)$ 
consisting of all smooth functions vanishing at a point $x\in U$, and 
$x=\pi (\ep\i(\ph\o i^*))=\pi(\ph(\pr_1\o i),\dots,\ph(\pr_n\o i))$, so 
that $\ep\i((i^*)^*(\ph))\in T_A(U)=\pi\i(U)\subset A^n$. 

Conversely for $a\in T_A(U)$ the homomorphism 
$\ep_a:C^\infty(\Bbb R^n,\Bbb R)\to A$ factors over 
$i^*:C^\infty(\Bbb R^n,\Bbb R)\to C^\infty(U,\Bbb R)$,
by steps 2 and 3 of \nmb!{4.5}. 

\smallskip
\noindent{\bf Step 3.} The two functors
$\operatorname{Hom}(C^\infty(\quad,\Bbb R),A)$ and $T_A:\Mf \to Set$ 
coincide on all open subsets of $\Bbb R^n$'s, so they have to
coincide on all manifolds, since smooth manifolds are exactly
the retracts of open subsets of $\Bbb R^n$'s see e.g\. \cit!{6},~1.14.1. 
Alternatively one
may check that the gluing process described in \nmb!{4.5},~step~4, 
%\cit!{6},~35.11,~step~6', 
works also for the functor
$\operatorname{Hom}(C^\infty(\quad,\Bbb R),A)$ and gives a
unique manifold structure on it which is compatible to $T_AM$.
\qed\enddemo

\enddocument